\documentclass{article}
\usepackage[nonatbib, final]{neurips_2019}

\usepackage{amssymb}
\usepackage{amsmath, amsthm}

\usepackage[utf8]{inputenc} 
\usepackage[T1]{fontenc}    
\usepackage{hyperref}       
\usepackage{url}            
\usepackage{booktabs}       
\usepackage{amsfonts}       
\usepackage{nicefrac}       
\usepackage{microtype}      

\usepackage{amsmath, amsthm}
\usepackage{setspace}
\usepackage{graphicx}
\usepackage{courier}
\usepackage{amsfonts}
\usepackage{amssymb}
\usepackage{algorithmicx}
\usepackage[noend]{algpseudocode}
\usepackage{algorithm}
\usepackage{wrapfig}
\usepackage{caption}

\newcommand{\bigA}{%
  \mathop{ 
    \mathchoice{\dobigA\Large}
               {\dobigA\large}
               {\dobigA\normalsize}
               {\dobigA\small}
    }\displaylimits 
}

\newcommand{\dobigA}[1]{%
  \vcenter{#1\kern.2ex\hbox{$\mathsf{A}$}\kern.2ex}}

\title{Analysis of heterogeneous computing approaches to simulating heat transfer in heterogeneous material}

\author{
 Andrew Loeb \\
Center for Applied Mathematics\\
Cornell University\\
Ithaca, NY 14850\\
  \texttt{ael89@cornell.edu} 
  \And
Christopher Earls\\
School of Civil and Environmental Engineering\\
Cornell University\\
Ithaca, NY 14850\\
  \texttt{earls@cornell.edu}
}

\begin{document}

\maketitle

\begin{abstract}
The simulation of heat flow through heterogeneous material is important for the design of structural and electronic components. Classical analytical solutions to the heat equation PDE are not known for many such domains, even those having simple geometries. The finite element method can provide approximations to a weak form continuum solution, with increasing accuracy as the number of degrees of freedom in the model increases. 
This comes at a cost of increased memory usage and computation time; even when taking advantage of sparse matrix techniques for the finite element system matrix. We summarize recent approaches in solving problems in structural mechanics and steady state heat conduction which do not require the explicit assembly of any system matrices, and adapt them to a method for solving the time-depended flow of heat. These approaches are highly parallelizable, and can be performed on graphical processing units (GPUs). Furthermore, they lend themselves to the simulation of heterogeneous material, with a minimum of added complexity. We present the mathematical framework of assembly-free FEM approaches, through which we summarize the benefits of GPU computation. We discuss our implementation using the OpenCL computing framework, and show how it is further adapted for use on multiple GPUs. We compare the performance of single and dual GPUs implementations of our method with previous GPU computing strategies from the literature and a CPU sparse matrix approach. The utility of the novel method is demonstrated through the solution of a real-world coefficient inverse problem that requires thousands of transient heat flow simulations, each of which involves solving a 1 million degree of freedom linear system over hundreds of time steps.
\end{abstract}


\section{Introduction}
\label{sec:introduction}

The  finite element method (FEM) provides useful approximations to weak solutions of partial differential equations (PDEs) over a much broader family of domains than those for which analytical solutions are known \cite{FEniCSBbook}. The method involves transforming a continuum solution into the solution of a sparse system of linear equations, which are well understood and can be quickly solved or approximated. This transformation encodes the domain geometry as well as properties of the domain that can vary over space and time, such as material properties for heat conduction problems. Fast solvers for this type of problem are critical for providing simulations over domains with increasingly fine discretization. This is particularly salient in the case that many successive simulations are desired with changes to the parameterization of material coefficients, for example: as in the solution of a coefficient inverse problem \cite{EarlsMSSP}. 

The present study is motivated by such inverse problems. Furthermore, increasingly efficient methods for the simulation of heat transfer through heterogeneous material are useful for the design of structural and electrical systems \cite{Kaviany}. In this area of research, analytical solutions have been found in some restricted cases \cite{Zoppou}, and only some characteristics of the transient behavior of more general problems are known \cite{Jaeger, Ash,Tsao}. The purely mathematical analysis relies on isotropy or bulk approximation over the domain. The presence of more than one set of heat conduction parameters over a heterogenous material makes the problem too complex to solve, or results in solutions that are too complicated to be useful \cite{Jaeger}.
 
\subsection{Scope and Organization}
\label{subsec:scopeAndOrganization}

This paper is divided into six sections. In Section \ref{sec:backgroundAndMotivation}, some history of matrix-free methods for FEM is discussed; leading to recent heterogeneous computing approaches. In Section \ref{sec:problemDescription}, the necessary mathematical preliminaries for our methods are summarized. The assembly operator is introduced, along with the flexibility that it offers for matrix-free methods. This flexibility is explored in Section \ref{sec:Implementation} with three interpretations of the assembly operator and the resulting implementations for the simulation of heat transfer through heterogeneous media. In Section \ref{sec:resultsAndDiscussion}, numerical experiments and their results are described, to compare the speed of the implementations over a range of problem sizes. A sparse serial method is also included in the comparison. A coefficient inverse problem arising from a real-world need for nondestructive corrosion detection is introduced and analyzed in Section \ref{sec:MCMC}. Finally, Section \ref{sec:conclusions} summarizes the goals and future impact of this research.

\section{Background and Motivation}
\label{sec:backgroundAndMotivation}
The methods that are developed in this work are derived from an element-by-element (EbE) decomposition of the finite element method. This viewpoint was introduced by Hughes et al. in 1983, for heat conduction calculations that were large for the time \cite{Hughes_JEM}. The motivation at that time was in avoiding exceedance of available computer memory for the storage of the large FEM system matrix, rather than considerations of speed. The authors showed stability for the algorithm, and followed up to show the same properties of the EbE method for structural and solid mechanics problems \cite{Hughes_CMAME}. Carey et al. exposed the potential for the EbE method to be parallelized, with a demonstration of a 2D convection-diffusion simulation in 1988 \cite{Carey}. Following this, research interest in the EbE framework had little room for advancement until advances in computer hardware were made.  

Heterogeneous computing is the practice of using dissimilar coprocessors in the solution of a numerical problem; typically a computer's central processing unit (CPU) and one or more graphics processing units (GPUs). This provides the user with access to fast serial processing power, as well as large-scale parallelism for certain computations that are properly amenable. Kiss et al. first utilized GPUs for the EbE method, using NVIDIA's proprietary CUDA platform \cite{Kiss}. The authors explored considerations that are particular to GPU computing, such as the preference for repeated computation rather than loading data from memory, and segmenting the domain with graph-coloring methods to avoid race conditions in parallelism. There have since been other efforts to parallelize PDE solvers for GPUs that are leveraged on symmetries of a particular problem \cite{Muller}. A noteable recent package is the TeaLeaf mini-app, which solves the heat equation in 2D using a matrix-free finite difference discretization \cite{McIntoshSmith2017}.

Modifications to the classical EbE decomposition have recently been explored by Mart\'inez-Frutos et al. \cite{Martinez-Frutos_CS}. The authors took a finer-grained approach to consider each degree of freedom (DoF) rather than a complete element, to formulate the DoF-by-DoF (DbD) method. It was shown in the context of elasticity problems that  synchronization overhead from graph coloring is more costly than the unfavorable memory access patterns which it prevents, especially for 3D problems. Mart\'inez-Frutos and Herrero-P\'erez further explored the DbD method for elasticity problems on a fixed grid mesh so that only one local finite element matrix is required  \cite{Martinez-Frutos_FEAD}. By varying material coefficients at the element level, between a constant ``inside'' and zero ``outside'', different domain geometries were enforced on the same mesh. More recent work by the same authors applied these strategies to problems of robust topology optimization \cite{Martinez-Frutos_CMAME}. By leveraging the ability to simulate many different domain geometries, an optimal  structural design for a given set of loading conditions can be found. Heterogeneous computing methods for topology optimization have also been studied for problems of designing domains with desired thermal properties. These studies have shown good results for problems involving steady state heat transfer in 2D \cite{Wadbro2009} and 3D \cite{Martinez-Frutos_CMAME}, using FEM solutions to the elliptic time-independent heat equation. Lastly, Mart\'inez-Frutos and Herrero-P\'erez have shown that multiple GPUs can be used effectively in the solution of topology optimization problems through task-level parallelism---the simultaneous evaluation of independent models that arise within a collocation strategy. All of the heterogeneous computing research described above, with the exception of TeaLeaf, was done with the CUDA platform.

The extension of heterogeneous computing beyond a single GPU has been an active research area, as further advances in scientific computing will require these paths. For example, Yang et al. study a general approach to partitioning the task of sparse matrix-vector multiplication across a range of heterogeneous architectures \cite{Yang2015}. A strength of the approach is in emphasizing the mathematical model separate from the computational architecture, to allow greater flexibility. Similarly, Gao et al. optimize a matrix-explicit conjugate gradient solver for multiple GPU, with attention given to five candidate sparse matrix-vector multiplication formats to automatically select the one that is most appropriate for the hardware \cite{Gao2017}. Considering that this is, by far, the most computationally expensive stage of the algorithm, the attention to the choice of memory storage format is warranted.

The present paper describes the theory and implementation of three approaches to a matrix-free preconditioned conjugate gradient algorithm for simulating transient heat conduction through a heterogeneous medium. Two are guided from previous studies, while our third combines benefits from both. The implementations differ in the interpretations of the DbD decomposition, as well as varying the use of a fixed grid or a general mesh, and in considering specialized hardware capabilities of GPUs. The most advanced implementation uses coalesced transactions with global memory for  improvements in hardware efficiency, and is modified through a domain decomposition to run across dual GPUs. The domain decomposition is done in a way that minimizes communication between the two devices, so that the additional computational power can be effectively deployed. All of the implementation are made within the OpenCL computing framework, which is non-proprietary and free to use on any platform \cite{OpenCL}. Scripting is done with the PyOpenCL package, providing readability and convenience with virtually no sacrifice in performance of the OpenCL API \cite{PyOpenCL}. The performant code is available for public use, distribution, and modification \cite{LoebGit}. 

\section{Problem description}
\label{sec:problemDescription}
We present the mathematical and computational context for assembly-free finite element methods with a focus on the parabolic time-dependent heat equation PDE.
\subsection{FE Formulation I (PDE)}
We wish to solve the heat equation in three dimensions with spatially dependent material coefficients. In strong form, the boundary value problem is
\[\begin{cases}
\rho C\frac{\partial T(\vec{x}, t)}{\partial t} = \nabla \cdot(k \nabla T(\vec{x}, t)) & \mbox{ in domain } \Omega \times (0, t_f),\\
k \frac{\partial T(\vec{x}, t)}{\partial \vec{n}} = f(\vec{x}) & \mbox{ on sides,} \\
 T(\vec{x},0) = T_{\text{ambient}},
\end{cases}\]
where the relevant thermal properties are the material {density}, $\rho$, {specific heat}, $C$, and {thermal conductivity}, $k$, all of which are assumed to be constant with respect to temperature. Discretizing in time with a $\theta$-scheme \cite{FEniCSBbook}, the spatio-temporal temperature profile $T(\vec{x}, t)$ is reduced to a finite set of temperatures at regularly spaced time increments, $\{T^{(i)}(\vec{x})\}_{i\in \mathcal{I}}$, $\mathcal{I} = \{0, 1, \dots, t_f/\Delta t\}$. The problem is then converted to weak form by multiplying the strong form with a test function $\phi(\vec{x})$ and integrating by parts to give the operators
\begin{align*}
a(T^{(i)}(\vec{x}),\phi(\vec{x})) &= \int_\Omega \left(  \rho C T^{(i)}\phi + \theta \Delta t k \nabla T^{(i)} \cdot \nabla \phi\right) d\vec{x} \\
L(\phi(\vec{x})) &= \int_\Omega \left( \rho C T^{(i-1)}\phi - (1- \theta) \Delta t k \nabla T^{(i-1)} \cdot \nabla \phi \right) d\vec{x} + \int_{\partial \Omega} \Delta t  f \phi ds.
\end{align*}
We presume $T^{(i-1)}$ to be known, and find $T^{(i)}$ so that $a(T^{(i)},\phi)=L(\phi)$ for all $\phi$ in some family of functions. If this family is composed of a finite set of basis functions $\{\phi_m\}_{m\in \{1,\dots,N\}}$, the finite element method can be used to solve for $T^{(i)}$ as a linear combination of them: $T^{(i)} = \sum_{m} U^{(i)}_m \phi_m$, where $\vec{U}$ is a vector of coefficients. The discrete boundary integral of $f$ is computed to give the vector $\vec{F}$. The finite element method involves the resulting matrices of pairwise integrals of basis functions
\[ \mathbf{M} = \left[ \int_\Omega \rho C \phi_m \hat{\phi}_n d\vec{x}  \right]_{m,n\in \{1,\dots,N\}} \text{ and } \mathbf{K} = \left[ \int_\Omega k \nabla  \phi_m \cdot \nabla \hat{\phi}_n d\vec{x} \right]_{m,n\in \{1,\dots,N\}}.  \]
Computation of these matrices is the process of \textit{assembly}. With $\mathbf{M}$ and $\mathbf{K}$ available, solving the weak form of the heat equation for all $\phi \in \{\phi_m\}_{m\in \{1,\dots,N\}}$ is equivalent to solving the matrix equation at each time step
\[ \left[ \mathbf{M} + \theta \Delta t \mathbf{K} \right] \vec{U}^{(i)} = \left[ \mathbf{M} - (1-\theta)\Delta t  \mathbf{K} \right] \vec{U}^{(i-1)} + \Delta t \vec{F} \]
for $\vec{U}^{(i)}$. This differs from the process of solving for a steady state heat flow in two ways. First, the FEM system matrix has a more complicated structure, rather than only involving the stiffness matrix $\mathbf{K}$. Second, the system must be solved at every time step to produce a transient solution. The time discretization process requires its own considerations for numerical accuracy and stability \cite{FEniCSBbook}. In this work we set $\theta = 0.5$, corresponding to a Crank-Nicolson method.

\subsection{FE Formulation II (Assembly-Free Methods)}
The matrices $\mathbf{M}$ and $\mathbf{K}$ are typically constructed by summing local contributions from each element in the assembly process. A local assembly matrix for element $e$, with $D$ degrees of freedom, contains the pairwise inner products of all basis functions with support in element $e$,
\[ \mathbf{M}_e  = \left[ \int_{\Omega_e} \phi_m \hat{\phi}_n d\vec{x}  \right]_{m,n\in \{1,\dots,D\}}. \]
Material property coefficients are taken to be constant over each element, so that the elemental assembly matrices depend only on the geometry of the domain. If all of the finite elements are the same size and shape, a single elemental assembly matrix can be reused and the mesh is said to have a \textit{fixed grid} (FG). The \textit{assembly operator,} $\bigA$, over the index set of elements $\mathcal{E}$ denotes the process of constructing a full system matrix from its local contributions. For example, 
\[ \bigA_{e\in \mathcal{E}} (\rho C)_e \mathbf{M}_e = \mathbf{M}. \] 
The assembly operator can also be applied to contributions within a single vector over each element to give the full vector. It is only a notational convenience to describe the mapping from local degrees of freedom to sums over global degrees of freedom. As such, the following are  valid notation for the general expression $\mathbf{M}\vec{x} = \vec{y}$
\begin{equation}
\label{eq:EbEvsDbD}
\bigA_{e \in \mathcal{E}} (\rho C)_e \mathbf{M}_e \vec{x}_e =  \bigA_{n \in \mathcal{N}}\left(\sum_{e\in \mathcal{E}^{(n)}} (\rho C)_e \mathbf{M}_e^n \vec{x}_e \right) = \vec{y}, 
\end{equation}
meaning that assembly is computed in terms of the degrees of freedom (over index set $\mathcal{N}$) in the ``outer loop'' with each of their elemental contributions computed separately. The first method is an EbE approach, similar to the standard method of assembling $\mathbf{M}$. The second is a DbD approach, in which the necessary vector dot products are viewed with finer granularity \cite{Martinez-Frutos_CS}. The freedom of interpretation of the assembly operator gives rise to the different strategies for parallel matrix-vector multiplication that are described in Section \ref{sec:Implementation}. 

To simplify notation, let $\mathbf{A} = \left[ \mathbf{M} + \theta \Delta t \mathbf{K} \right]$ and $\mathbf{L} = \left[ \mathbf{M} - (1 - \theta) \Delta t \mathbf{K} \right]$. Set $\vec{b} = \mathbf{L}\vec{U}^{i-1} + \vec{F}$ and $\mathbf{A}_e = (\rho C)_e \mathbf{M}_e + \theta \Delta t k_e \mathbf{K}_e$ for $e \in \mathcal{E}$. Then the problem of finding $\vec{U}^{i}$ at each time step is reduced to solving
\[ \mathbf{A} \vec{U}^i = \bigA_{e \in \mathcal{E}} \mathbf{A}_e \vec{U}^{i}_e = \bigA_{n \in \mathcal{N}}\left(\sum_{e\in \mathcal{E}^{(n)}} \mathbf{A}_e^n \vec{U}_e^{i} \right) = \vec{b}.\]
We note once again that the explicit computation and storage of $\mathbf{A}$ and $\mathbf{L}$ is not necessary if $\{\mathbf{M}_e\}_{e\in \mathcal{E}}$ and $\{\mathbf{K}_e\}_{e\in \mathcal{E}}$ are available. 

\begin{figure}
\centering
\includegraphics[width=\textwidth]{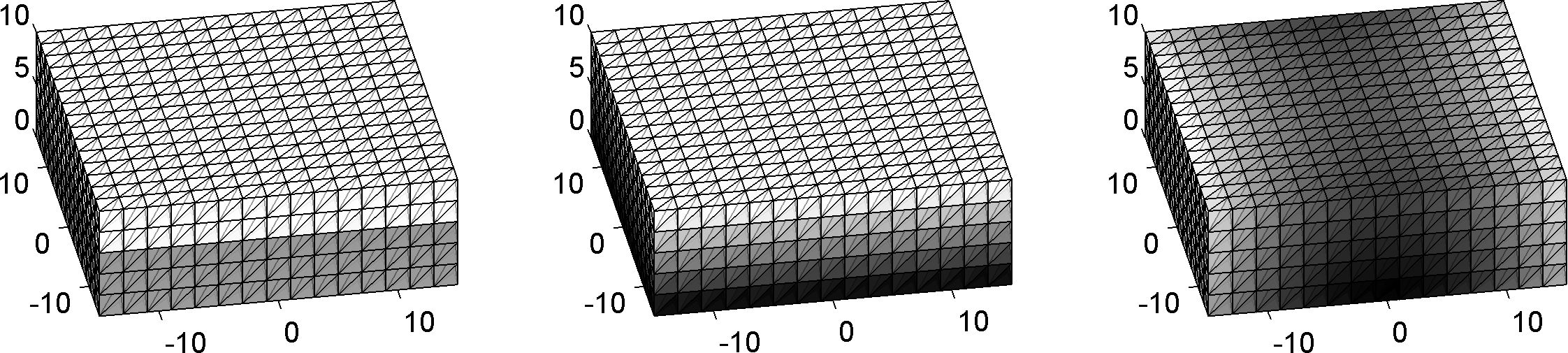}
\captionsetup{margin = 0.5cm}
\caption{30$\times$30$\times$10 mm domain with tetrahedral meshing and three functions describing changes in material properties. A boundary between two materials can be abrupt (left) or smoothed over many elements (middle). Heat conduction can also be simulated on a domain with more complex dependence between material properties and spatial location (right) with negligible computational burden, provided that the properties have a closed functional form (here, the shading is computed as $x^2 - 0.2y^2 + 10z$ at each vertex and the values averaged over each element).}
\vspace{-.4cm}
\label{fig:domain}
\end{figure}

Assembly-free methods are especially useful if the spatially dependent material properties are not known in advance, or if many simulations are to be done over the same domain with varying coefficients. The generation of the mesh geometry and the computation of elemental assembly matrices can be done in advance and stored. Then all of the remaining computations required for a matrix-vector multiplication are parallelizable. Figure \ref{fig:domain} illustrates a 3D FG mesh with tetrahedral elements and three sets of material properties parameterized by shading. Each of these spatially varying functions for the material properties have a simple functional form, and can be efficiently employed with elementary assembly matrices that are computed and stored beforehand,  to treat changing analysis contexts under material property variation.

\subsection{Preconditioned Conjugate Gradient}
\label{PCG}
The matrix $\mathbf{A} = \left[ \mathbf{M} + \theta \Delta t \mathbf{K} \right]$ is large, sparse, symmetric, and positive definite. The system above is thus solvable with the preconditioned conjugate gradient (PCG) algorithm \cite{Shewchuk}. At each time step, let $\vec{x} = \vec{U}^{(i)}$ (to free superscripts and the index $i$), be set at some initial guess. There also must be provided $\mathbf{A}, \vec{b}$, a residual tolerance, and a preconditioner matrix $\mathbf{P}$ for which $\mathbf{P}^{-1} \vec{x}$ is easily computable and $\mathbf{P}^{-1}\mathbf{A}$ is relatively well conditioned, such as the Jacobi preconditioner. The algorithm involves the following linear algebra operations: matrix-vector multiplication (MVM), diagonal inverse matrix-vector multiplication (DIMVM), vector-vector multiplication (VVM), and adding scalar multiples of vectors (AXPY). The dominating computation is the MVM in each iteration. All of the other computations are easily parallelizable. The parallelization and implementation of the MVM on GPUs for system matrices of heat conduction FEM problems is the focus of this work. The subsequent PCG algorithm is shown in Algorithm \ref{alg:PCG} with each step annotated by its type of linear algebra operation. 
\begin{algorithm}
  \caption{Preconditioned Conjugate Gradient}
\label{alg:PCG}
  \begin{algorithmic}[1]
    \Statex
    \Function{PCG}{$A, b, x, i_{\max}, \text{tol}, P$}
\State  $i \gets 0$
\State $r \gets b - Ax$                              \Comment{MVM}
\State $d \gets P^{-1}r$                               \Comment{DIMVM}
\State $\delta_{\text{new}} \gets r^Td $                       \Comment{VVM}
\While {$(i < i_{\max})$ and $(\delta_{\text{new}} > \text{tol})$ }
   \State $q \gets Ad$                               \Comment{MVM}
  \State  $\alpha \gets \delta_{\text{new}}/(d^Tq)$              \Comment{VVM}
  \State  $x \gets x + \alpha d$                       \Comment{AXPY}
   \If {$i$ is divisible by 50}
     \State $r \gets b - Ax$                        \Comment{MVM}	      
   \Else
      \State $r \gets r - \alpha q$                    \Comment{AXPY}
\EndIf
   \State$s \gets P^{-1}r$                            \Comment{DIMVM}
   \State$\delta_{\text{new}} \gets r^Ts $                     \Comment{VVM}
   \State $\beta \gets \delta_{\text{new}}/\delta_{\text{old}}$
   \State$d \gets s + \beta d$                       \Comment{AXPY}
   \State$i \gets i + 1$
   \EndWhile
    \EndFunction
  \end{algorithmic}
\end{algorithm}

A Jacobi preconditioner is used for all of the methods described in this work \cite{Shewchuk}. The construction of a Jacobi preconditioner is a straightforward process that lends itself to element-wise parallel computation. We note that the development of preconditioners that can be computed on a GPU is an active area of research \cite{Martinez-Frutos_FEAD, Golub1999}, but is not a focus of the present study. A comparison of this strategy with a serial implementation having a stronger preconditioner is made in Section \ref{sec:serialComparison}.

\subsection{OpenCL Heterogeneous Computing Framework}
Here, we summarize the framework of OpenCL computing to introduce common language that will be used in Section \ref{sec:Implementation}.
\subsubsection{Computation hierarchy}
With OpenCL, GPUs are programmed with \textit{kernels}; small bits of C code that are sent in parallel to the individual cores \cite{OpenCL}. At any given time, each of the many cores is acting as a \textit{work item}, which is the most granular operating unit in the hierarchy. Work items are structured in \textit{work groups}, with as few as one work item per work group. The user is responsible for defining the sizes and dimensions of the hierarchical structure, so that at the time of execution, each work item is provided with unique identifying information and the generic kernel code. The identifying information is:
\begin{itemize}
\item \textit{global\_id} ranging from 0 to the total number of work items in each dimension,
\item \textit{local\_id} ranging from 0 to the number of work items in a work group, in each dimension,
\item \textit{group\_id} ranging from 0 to the total number of work groups in each dimension.
\end{itemize}
The total number of work items, total number of work groups, and size of each work group is also available.  The kernel explicitly tells each work item how to contextualize itself within the larger problem, to determine what data must be loaded from memory or computed privately. At the end of the kernel, results of the independent granular computations are written to memory.

\subsubsection{Memory hierarchy}
\label{sec:MemoryHierarchy}
The use of memory on a GPU dictates programming strategies and the success or failure of an algorithm. At fully efficient throughput, a single core can execute several floating point operations per clock cycle \cite{OpenCL_bestPractices, du2012cuda}. However, accessing data from the ``slow'' global memory location can take 400-600 clock cycles. This alone warrants special attention to the OpenCL memory hierarchy.

Data that is loaded onto the GPU, or stored as the output of work items, must be stored in \textit{global} memory, which has space on the order of gigabytes. During execution, all work items have their own small amount of \textit{private} memory, which is on-chip and not visible to any other work item. This is fast to access, and used for variables that can take different values across every work item. In between, there is \textit{local} memory, which is also on-chip, and shared among a single work group. There are usually tens of kilobytes reserved for local memory for each work item. Access to local memory is roughly 100 times faster than accessing global memory, provided that work items within the work group aren't making conflicting calls (``bank conflict''). Local memory is allocated outside of the kernel, and cannot be freely initialized with specified data. Finally, kernels can consider certain data as \textit{constant} memory. This data is physically still in global memory, but a kernel cannot write to it. When a kernel reads data from constant memory, it is cached, so that subsequent reads are fast. Global memory is not cached.

There are two strategies for gracefully managing reads from global memory when it is necessary. First, the latency can effectively be hidden if there is enough non-dependent computation to keep a work item busy between the time when the data are called and the time they are used. This is preferable, though not always possible. Second, a kernel can take advantage of the way that hardware loads data from global memory to private memory. A single transaction with global memory yields 32 words (such as 8 byte double-precision floats) of data, whether it is all called for or not. If work items that are indexed sequentially by global ID, request data from global memory that is organized in the same sequential way, the calls are automatically bundled and processed as one transaction. This process is the simplest form of \textit{coalesced} memory access. There have been advances in hardware, and in OpenCL standards, to provide more flexibility, such as allowing permutations of the 32 sequential words to 32 work items, that are blocked together but not necessarily in the same order.

\subsubsection{Programming strategy}
The OpenCL programming approach is as follows:
\begin{enumerate}
\item Investigate the hardware to find and define a \textit{host} (CPU/hard drive etc., the conventional computing environment), its \textit{devices(s)} (GPU with its onboard memory), and define a \textit{context} the overall computing environment. 
\item Define initial variables on the host.
\item Load data onto a device. This includes reserving space for data that a kernel will write later and anything that is meant to last from one kernel to another. The benefit of defining all of the memory space in advance is that the user may specify whether each buffer is read or write only (or both) for both the host and device, or even if it is known that the host will never try to read it. Then the space that is allocated will be optimal for however the data will be treated.
\item \textit{Build} the compute kernels. This involves compiling the C code and specifying pointers to memory buffers where its arguments can be found. A single kernel program can be built multiple times with different arguments, as is the case with the vector-vector multiplications in steps 5, 8, and 15 in Algorithm \ref{alg:PCG}. Each of these are built independently.
\item Define a \textit{queue} in the context. The host can enqueue kernels, memory transfer operations, or wait fences. OpenCL turns these into individual tasks that are performed as cores on the GPU become available. Code is written as if the context has infinitely many cores to run in parallel, and then the queue manages the execution of code on available hardware. Flags can be used with enqueued commands to ensure that all tasks from one kernel are finished before any tasks from the next kernel start, in case memory is being written and then read in a dependent way.
\item Enqueue commands to copy memory from a device to the host. This can be the final result of computations, or in the case of this work, the PCG residual, so that the host can decide whether or not to begin another iteration of enqueueing kernel commands.
\end{enumerate}

\section{Implementation of Assembly-Free Methods}
\label{sec:Implementation}
We outline three assembly-free algorithms for matrix-vector multiplications. The differences arise from the flexibility in interpreting the assembly operator demonstrated in Equation (\ref{eq:EbEvsDbD}). For each method, the explicit assembly equation is provided, along with an outline of the memory and computational hierarchy for a parallel implementation on GPUs. Further details are discussed in \ref{app:PCG} so that broader concepts behind the implementations can be the present focus. We begin by giving context of the particular geometry of the problem.

\subsection{Mesh Geometry}
A 3D regular mesh with linear tetrahedral elements is generated based on the domain boundaries and the number of divisions in each dimension. The domain is then divided into rectangular prisms according to these divisions. The examples here are all cubes for simplicity. Each cube is then subdivided into six tetrahedra, as shown in Figure \ref{fig:domain}. A benefit of using tetrahedral elements is that their elemental assembly matrices have size $4\times 4$.  The OpenCL specification allows dot products of  the 4 element floating point vector data type, \textit{double4}, with a single instruction. This is the central operation for all of our fast MVM methods, regardless of the interpretation of the assembly operator.

\begin{figure}
\includegraphics[width=\textwidth]{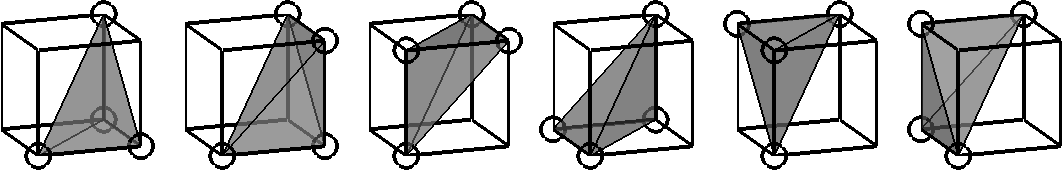}
\captionsetup{margin = 0.5cm}
\caption{Emphasized view of the tetrahedral subdivisions of each small cube in the mesh. Each finite element is characterized by four vertices of the cube.}
\vspace{-.4cm}
\label{fig:DoFMapLocal}
\end{figure}

The six tetrahedra within each cube are indexed in a sequential way---six in the first cube, then six in the next cube in the $x$ direction, etc. until the $y$ dimension is incremented, and then the next slice in the $z$ dimension begins after that. An emphasized view of the six tetrahedra is given in Figure \ref{fig:DoFMapLocal}. Each cube has consistent indexing as to which of its eight corners correspond to each element within; allowing for a single reference table to be stored in memory. The structure of this kind of mesh also permits local computation, such as determining spatially dependent material properties, without the need for a table of $xyz$-locations of each vertex. This information can be determined from the one-dimensional index of the vertex along with a small table of the domain boundaries and the number of divisions in each dimension. Substituting small computation from general domain information, in place of loading the same data from a lookup table, is important for parallelizing an algorithm for a GPU.

\subsection{Previous Strategies}
We begin by providing the details of two approaches that use ideas from previous literature.

\subsubsection{Implementation 1: Flexible DbD}
\label{sec:EbE}
The first implementation follows a flexible DbD strategy that does not assume that the FEM mesh lies on a fixed grid, but one which can be generated by a nonlinear transformation of one. Each global DoF of the vector $\vec{y} = \mathbf{A}\vec{x}$ in a general matrix equation is computed as
\begin{equation}
\label{eq:EbE}
\vec{y}_i = \left[ \sum_{\substack{e \in \mathcal{E}^{(i)} \\ i = \mathcal{N}^{(e)}(j) }} \left[ \mathbf{A}_e^j \vec{x}_e \right] \right] \; (i \in \mathcal{N}),
\end{equation}

\begin{wrapfigure}[22]{r}{0.5\textwidth}
\centering
\centerline{\includegraphics[width=0.45\textwidth]{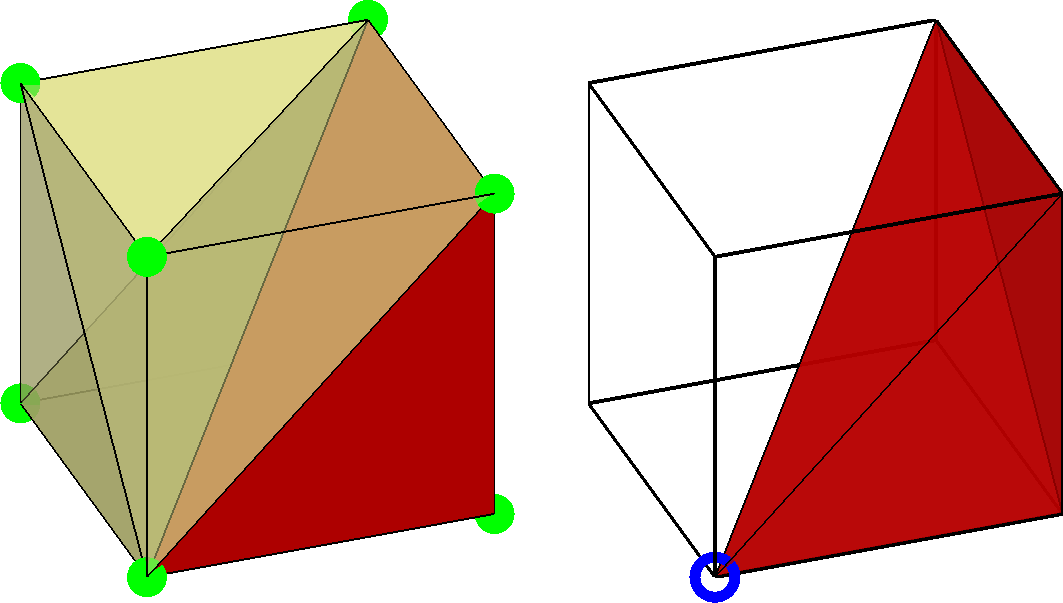}}
\captionsetup{margin = 0.1cm}
\caption{(left) Local data requirements for the first assembly-free MVM method. Elemental assembly matrices are required for six tetrahedral elements that comprise a cube in the mesh, as well as entries of the input vector corresponding to its eight corners. One element is emphasized and isolated (right), denoting the responsibilities of a single work item. There are 24 work items, within each work group for all such finite element-DoF pairs, required to assemble the output vector.}
\label{fig:EbEMem}
\end{wrapfigure}

where each pair of square brackets denotes a computation done by one work item in the implementation. Every elemental assembly matrix is stored in global memory, with a preprocessing step that determines $\{\mathbf{A}_e\}_{e\in \mathcal{E}}$  and $\{\mathbf{L}_e\}_{e\in \mathcal{E}}$ based on local material properties. A subset of this computation furnishes the Jacobi preconditioner matrix at the same time, by only storing diagonal entries. Work groups are responsible for each set of six elements in a cube, with 24 work items per work group, each corresponding to one element-DoF pair. Elemental assembly matrix coefficients are computed one time, and the scaled matrices are stored in double4 vectors in ``element order'' within global memory. In the computation, each of the 24 work items loads vertex data to local memory, including duplicated vertices, to alleviate conflicting simultaneous memory access from multiple work items. Then it reads its row from the elemental assembly matrix data, takes necessary entries from local vector data, performs a double4 dot product, and stores the result as its contribution to the global vertex in ``vertex order''. In a second pass, one work group for each global vertex reads these 24 consecutive contribution and sums them into the single entry of the result vector. A similar two step approach is described by Mart\'inez-Frutos et al. as an alternative to atomic incrementation to a single memory location for each degree of freedom \cite{Martinez-Frutos_CS}.

The necessary data and computational responsibilities of the first pass of this process are illustrated in Figure \ref{fig:EbEMem}. The left panel shows the scope of each work group---six elemental assembly matrices and eight entries of the input vector $\vec{x}$, which are loaded into the work group's local memory. The right panel illustrates the responsibilities of each work item---one element-DoF pair, corresponding to a length-four vector dot product that it must compute and store. For clarity and consistency with the following sections, the necessary data for this single work item is emphasized in the left panel as well.
 
Splitting the assembly operator into two explicit sums requires a sort of data transpose at some point of the computation. This is because assembly matrix data is naturally stored in element order, while the output must be combined and stored in vertex order. The two step approach here ensures efficient data retrieved from global memory for both steps. The consequences of writing data to potentially distant locations at the end of the first step are hidden, since the data is enqueued to be written, and then the kernel can be restarted and the computation continues on while the writing takes place. Furthermore, we note the second step can be slightly modified to perform an extra vector-vector addition of the form $ \vec{y} = \mathbf{A}\vec{x} + \vec{b}$ in the same kernel, which is one of the operations of the PCG algorithm. 

\subsubsection{Implementation 2: Single Pass FG DbD}
\begin{figure}
\centering
\includegraphics[width=\textwidth]{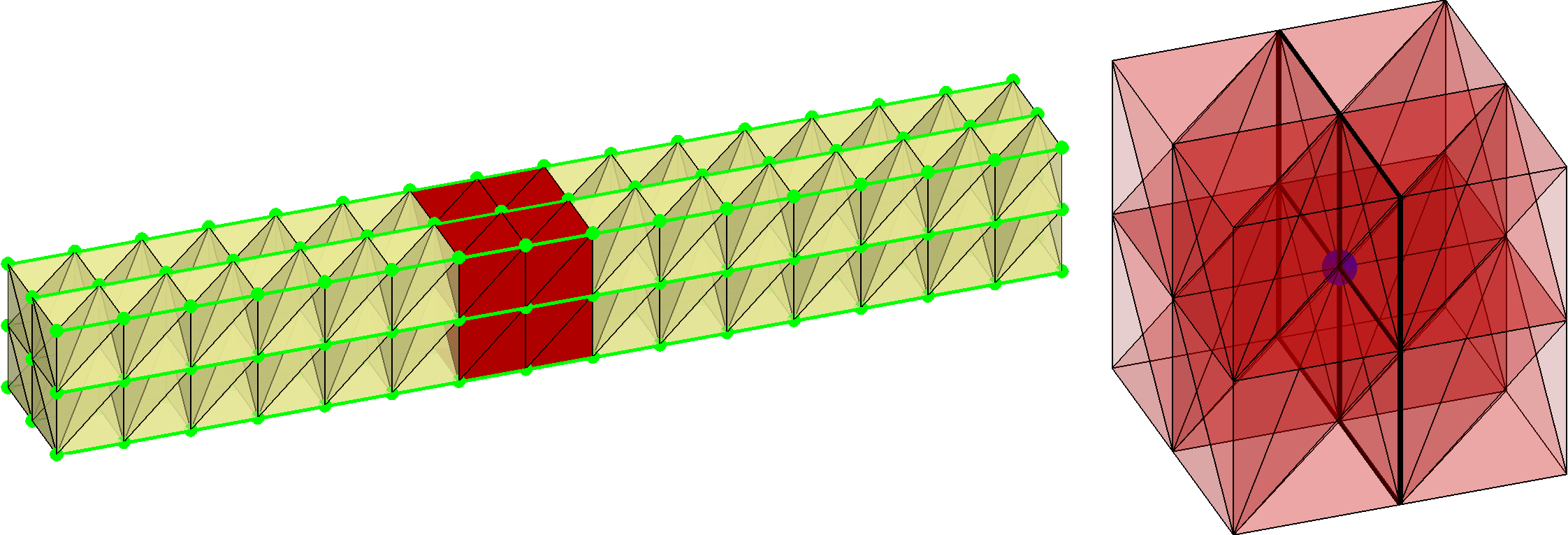}
\captionsetup{margin = 0.5cm}
\caption{Necessary local data and work item responsibilities for the second assembly-free MVM method. (left) Data that is loaded from adjacent locations in memory for the input vector are connected by a green line to emphasize the potential for coalesced memory loading. Note that this figure is truncated for clarity and that the method actually loads 64 consecutive entries to local memory. (right) A single work item requires more data, but fully computes an entry of the output vector, denoted by a filled point. As before, the data required by a single representative work item is emphasized from local memory on the left. }
\vspace{-.4cm}
\label{fig:DbDMem}
\end{figure}

The second implementation follows the coarser parsing of the assembly operator, with 
\begin{equation}
\label{eq:DbD}
\vec{y}_i = \left[ \sum_{\substack{e \in \mathcal{E}^{(i)} \\ i = \mathcal{N}^{(e)}(j) }}  \mathbf{A}_e^j \vec{x}_e  \right]  \; (i \in \mathcal{N}).
\end{equation}
Once again, the square brackets denote the computations of a single work item, so that this method only requires a single pass. Single pass, fixed grid DbD strategies have been previously explored in References \cite{Martinez-Frutos_FEAD, MartnezFrutos2017}. The trade off here is that this requires more data to be loaded into local memory for each work group for the full determination of a degree of freedom. Consequently, more memory must be loaded overall. This is partially alleviated by setting the work groups to be as large as allowed by hardware limitations, so that most of the data are reused by adjacent degrees of freedom. Previous strategies to maximize on-chip memory have organized the input data into 2D square patches, which can be loaded from global memory in a coalescent way \cite{Schmidt2011}. However, in order to have access to the data in neighboring elements, a halo of data around the patches must be loaded inefficiently. We choose to organize the necessary input data into 3 $\times$ 3 rectangular blocks, and as long as possible, so that all global memory access can be coalesced in the long direction, as shown in Figure \ref{fig:DbDMem}. For our hardware, specified in Section \ref{sec:resultsAndDiscussion}, this corresponds to 64 work items and 27 KB of local memory usage. We also implement this approach to be used on a fixed grid mesh so that only one elemental assembly matrix is required. This allows the elemental assembly matrix to be stored in constant memory, so that it does not need to be loaded anew by each work item. The Jacobi preconditioner is computed in the same way as was described in Section \ref{sec:EbE}, although the local material properties are not precomputed. This is not necessarily detrimental, since computing local material scaling coefficients is a small computation which can hide the latency of loading data from the input vector.

A visual overview of this method is provided in Figure \ref{fig:DbDMem}. Each work item is responsible for all contributions to one entry in the output vector, so it needs data from 24 local assembly matrices and 27 entries of the input vector. An internal loop scans through all 24 element-DoF contributions, computing elemental scaling coefficients, performing dot products, and cumulatively summing the result.

\subsection{Implementation 3: FG DbD with Memory Coalescing}
\label{sec:FGEbE} 
\begin{figure}
\centering
\includegraphics[width=\textwidth]{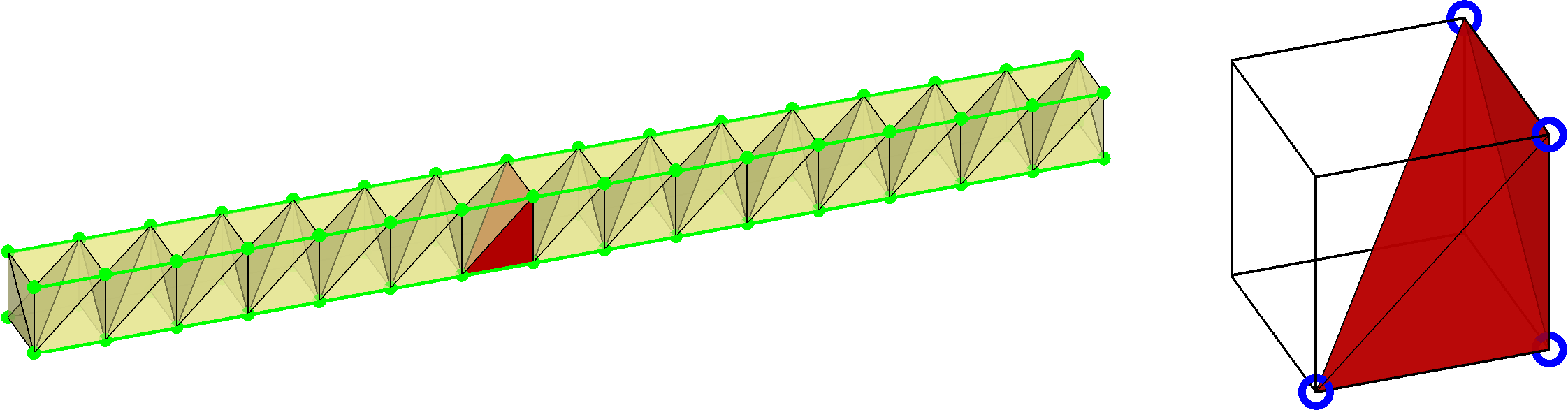}
\captionsetup{margin = 0.5cm}
\caption{Necessary local data and work item responsibilities for the third assembly-free MVM method. (left) Four coalesced reads from global memory provide all of the necessary input vector data. Note that this figure is truncated for clarity and that the method actually reads 32 consecutive entries to local memory with each coalesced memory read. (right) The representative work item computes contributions for all degrees of freedom associated with its finite element.}
\vspace{-.4cm}
\label{fig:EbEcMem}
\end{figure}

The third implementation combines advantages from the first two. It is similar in structure to the first, except with increased responsibility to each work item, larger work groups, and the restriction to a fixed grid. The explicit interpretation of the assembly operator for this method is the same as Equation (\ref{eq:EbE}), although the full elemental assembly matrix-elemental vector multiplication is carried out by a single work item here. The primary characteristic of this third method is that work groups are structured so that all data that is loaded from global memory in the first pass is done so with a coalesced access pattern, as described in Section \ref{sec:MemoryHierarchy}. In this sense, it is similar to the second implementation, except that the work group size is determined by the amount of data that can be loaded from a single coalesced memory read rather than by maximum local memory capacity. Furthermore, only four sections of memory are required for the first pass, as shown in Figure \ref{fig:EbEcMem}. The element-DoF contributions to each entry of the output vector are collected and summed as much as possible before writing back to global memory. While a second pass is still required to add the contributions among work groups, this process is faster than in the first implementation since there are fewer terms to sum. Additional details behind the efficient use of coalesced memory access now follow.

Since a fixed grid is assumed, the only data that must be loaded from global memory is the input vector. To do this efficiently, it must be loaded in blocks of 32 consecutive entries by blocks of 32 consecutive work items. A single memory read of this form does not give enough data to perform assembly for any element. However, if four blocks of memory are read, corresponding to the four horizontal edges of a long rectangular prism, then every tetrahedral element within that prism can be integrated over. That is, the contributions from each of these elements to the degrees of freedom that are loaded can be computed. 

\begin{figure}
\centering
\includegraphics[width=\textwidth]{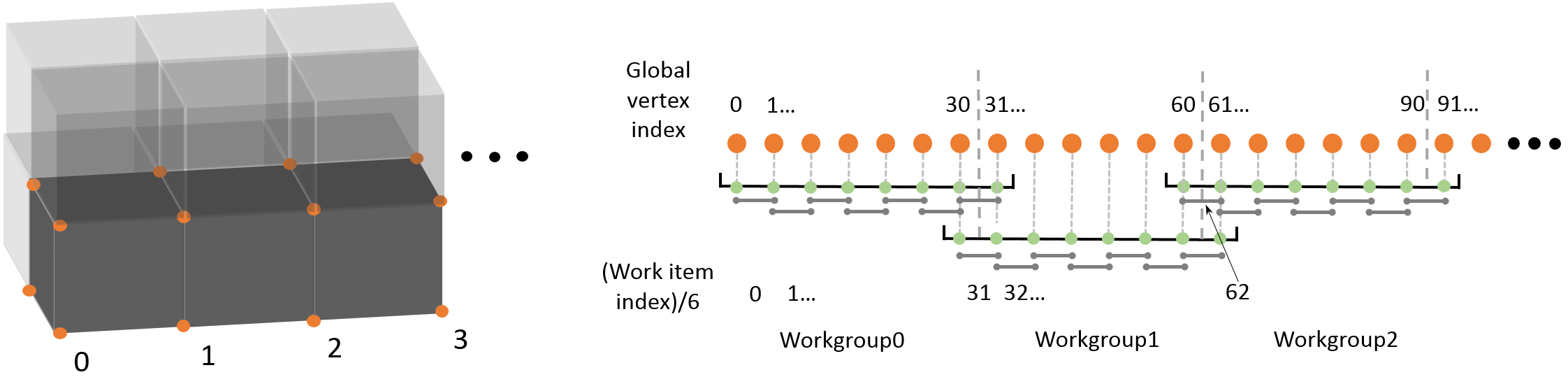}
\captionsetup{margin = 0.5cm}
\caption{Memory access and computational partitioning pattern for the coalesced DbD MVM method. The orange nodes correspond to data in global memory, while the green nodes represent local memory in each work group. The short gray lines within each work group represent the tetrahedral elements which must have access to the vertex data at all four corners to compute their assembly contributions. Four such coalesced reads from global memory are performed to provide the tetrahedra with their necessary data.}
\vspace{-.2cm}
\label{EbELoad}
\end{figure}

Using the memory access length of 32 as the guiding limit, 31 full cubes of six tetrahedral elements will be integrable. Therefore, work groups of 186 work items are invoked---one work item for each tetrahedral element. Since only 128 work items are necessary to read data from global memory, some work items are assigned a dual purpose and some remain idle during the loading process. Those first 128 receive both a global elemental index for their integration responsibility as well as a global vertex index for loading data from global memory into local memory. Figure \ref{EbELoad} depicts the way in which global data is accessed by each work item for one of the four edges of the long rectangular prism. The other edges are treated in the same way, with offset information determined in-kernel based on how many divisions are made in the domain in each dimension. The contribution to a vertex can only be computed if all four vertices of an adjacent element are included within a work group. Therefore, only the 30 internal vertices from each block of 32 is contributed to in storage (except for the first work group, and possibly the last). Vertices on either end comprise a halo region that must be loaded for the computation within the local memory geometry. This partitioning pattern is also demonstrated in Figure \ref{EbELoad}. The special treatment of halo data on the region boundary has been explored through several GPU computing methods, both the context of the global domain and for local computation \cite{McIntoshSmith2017, MartnezFrutos2017, Schmidt2011}. We next discuss our global memory halo.

\begin{figure}
\centering
\includegraphics[width=\textwidth]{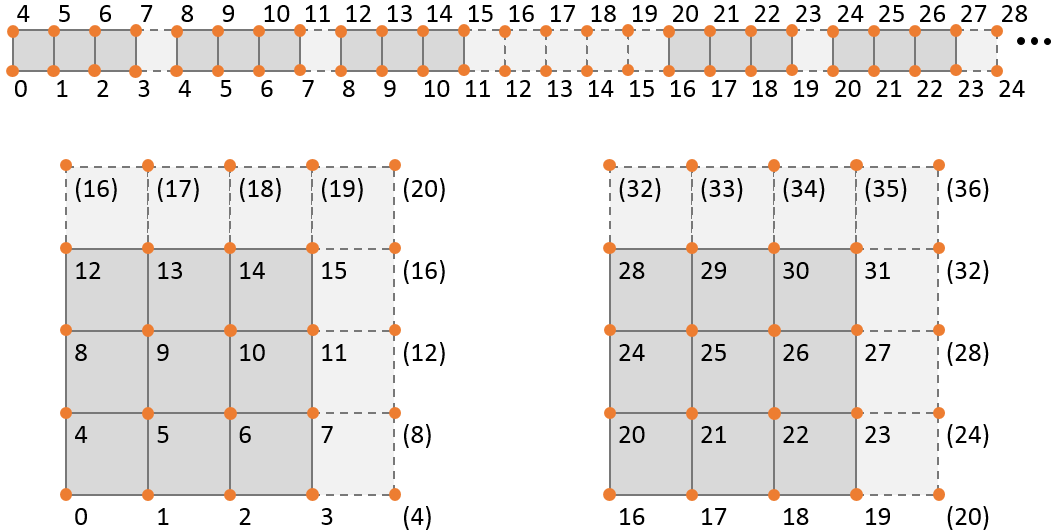}
\captionsetup{margin = 0.5cm}
\caption{Element padding for a 3$\times$3$\times$n mesh as seen from a top-down perspective. All of the data is loaded, but only contributions from solid elements are stored. This allows coalesced access to global memory to be used without interruption, at the cost of the extra discarded computation.}
\vspace{-.4cm}
\label{EbEStore}
\end{figure}

Another consequence of coalesced memory access is that some consecutive elements of the input vector do not actually share an element in the domain. This is demonstrated in Figure \ref{EbEStore} for a small example mesh. We avoid this problem by padding the domain with non-physical elements in the $+x$ and $+y$ directions from the perspectives of loading and computation, and then the contributions of these elements are ignored upon storage of the result. This produces some amount of wasted computation, but the relative volume of padding elements to the total mesh volume decreases as the granularity of the mesh increases. Furthermore, efficiency gains from coalesced memory access strongly outweigh the losses from this wasted computation. 

\subsection{FG DbD with Memory Coalescing on multiple GPUs}

\begin{figure}
\centering
\includegraphics[width=\textwidth]{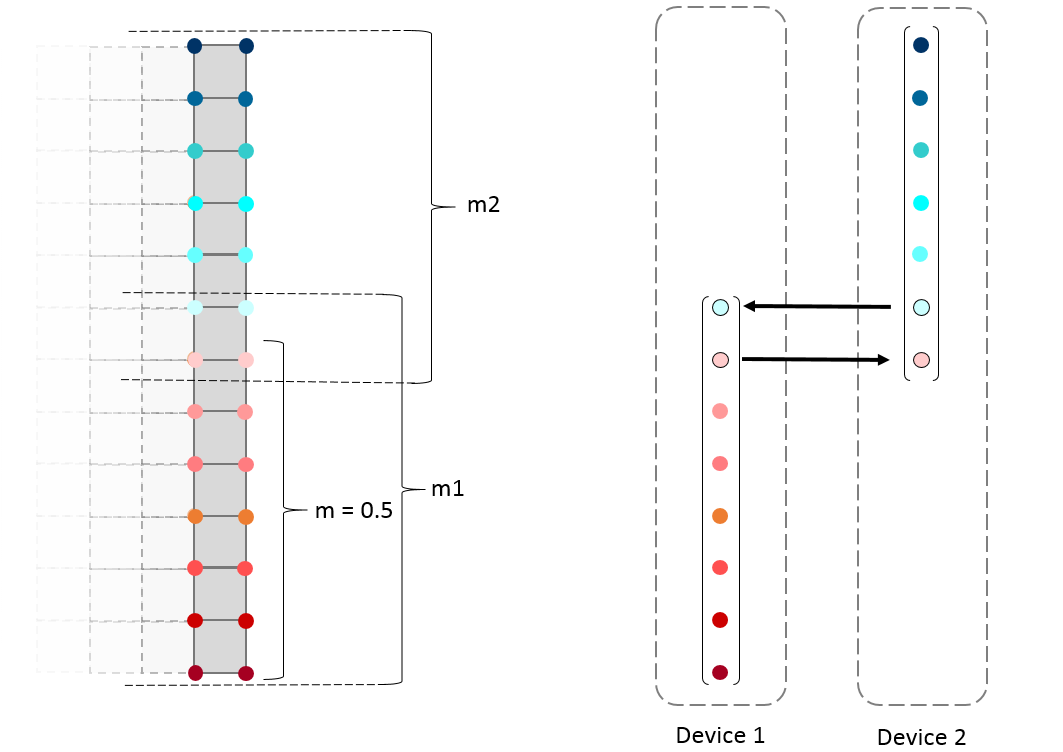}
\captionsetup{margin = 0.5cm}
\caption{(left) splitting of vertices between two devices according to a user-specified fraction $m$. (right) Data that must be transfered between devices at each iteration.}
\vspace{-.4cm}
\label{fig:dualGPU}
\end{figure}

The fixed grid coalesced method described above is modified for use on dual GPU. The domain is split in the $z$ direction according to the additive Schwartz method  \cite{Cai2003}. This keeps all vertices in each subdomain in adjacent blocks of memory. A fraction of the domain to be assigned to device 1, $m$, is specified, such as 0.5. This fraction of $z$ slices, rounded up, with one additional layer, is the number of vertices that device 1 receives for computing, $m_1$. The rest of the vertices, in addition to two overlapping layers are assigned to device 2, totaling $m_2$. This split is demonstrated in Figure \ref{fig:dualGPU}. If input vectors are initialized from a full set of global data, then one matrix-vector multiplication can be performed on each device, and the result can be faithfully reconstructed. For more than one sequential matrix-vector multiplication, the shared boundary data must be updated. Only one layer of vertices needs to be transferred in each direction. In the PCG algorithm, the solution vector $x$ is initialized at the beginning, and the intermediate vector $d$ must be transferred at each iteration. This memory transfer presents a bottleneck in the method, so that gains in speed are expected only with large systems for which MVM takes much longer than a GPU-to-GPU memory copy.

The scalar results of dot products must also be communicated between devices, both the residual $\delta$ and step size $\alpha$. To do this, each partial dot product is handled separately, omitting the extra boundary vertices. The partial results are stored in their own buffers, which are then transferred both directions. We thus require four memory buffers for each scalar quantity, one native buffer for the partial results on the device which computed it, and one target buffer on each device to store the copied value from the other. Many permutations of methods for the bidirectional communication of partial dot products were considered, such as copying the contributions to the host for summation. This method was found to be the fastest reliable way for each device to receive the full results. No special kernel is needed to combine them. Rather, modifications to the AXPY kernel are made to accept both buffers on a device and sum them as part of the existing process.

The only remaining matter for modification is adjusting the memory objects that contain vertex-to-spatial-location information. One each is made for the two devices, with the second encoding the $z$ shift for its first index. Then each device runs PCG in parallel, and can contextualize the location of vertices within the total domain. The host initialized two sets of every kernel that has been described for the serial method, as well as launching one queue to transfer boundary and scalar information and properly wait. Lastly, we note that the implementation here is for dual GPUs, but could be easily extended for increased distribution: with multiple GPUs. Additionally, the additive Schwartz method does not require uniformly sized decomposition regions. Therefore, heterogeneous GPUs could be used together, each taking a share of the total domain in order to balance the load across the devices.

\section{Experiments and Discussion}
\label{sec:resultsAndDiscussion}
The results shown here have been produced on AMD FirePro D700 GPU with 2048 streaming processors, 6GB of onboard memory, and up to 32KB of local memory per work group. The serial computations are done on a 2.7 GHz Intel Xeon E5 processor using a single core. The CPU algorithm is implemented in Python and optimized for solving heat conduction problems in many of the same ways as the GPU methods discussed previously. Elemental assembly matrices are precomputed so that assembly of the sparse system matrix can be done efficiently after the specification of material parameters, analogously to the parallel GPU algorithms. Furthermore, Jacobi preconditioning is specified so that the comparisons below are as fair as possible. 
 
\subsection{Performance Comparison}
\label{sec:performanceComparison}
\begin{figure}
\includegraphics[width=\textwidth]{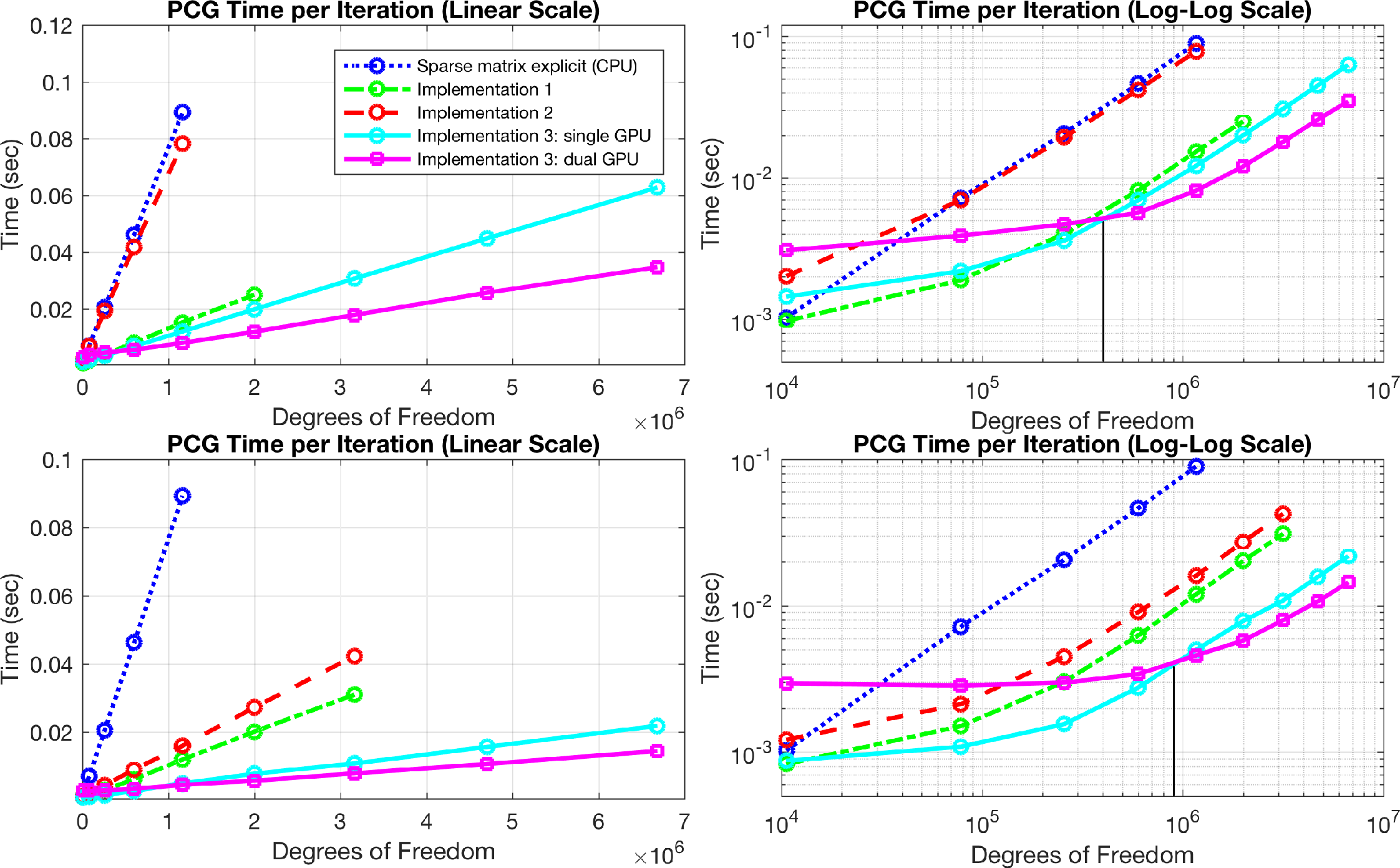}
\captionsetup{margin = 0.5cm}
\caption{(Top) Time per PCG iteration for four GPU implementations and a serial sparse implementation of the same algorithm. Results are computed for each method as far as hardware limitations would permit. (Bottom) Profiling results using GPU single-precision floating point operations. The serial sparse matrix results are the same as before for consistent comparison. The point at which the dual GPU implementation overtakes the single GPU implementation is marked.}
\label{fig:PCG_time_per_iteration}
\end{figure}

We first report the performance of each implementation over a range of problem sizes by simulating uniform heating on the front face of a two-layer laminate. The domain for this problem is a rectangular prism $\Omega = [-15, 15] \times [-15, 15] \times [0, 10]$ with $f(\vec{x}) = 1$ on $x_3 = 0$ and zero everywhere else, and $T_{\text{ambient}} = 0$. The material parameters are specified as
\[ \rho C = \begin{cases} 3.724\text{e}6 \text{ g/mm C s}^2 & x_3 \leq 5\\ 1.65\text{e}6 \text{ g/mm\,C\,s}^2 & x_3 > 5 \end{cases},\qquad k = \begin{cases} 4.9\text{e}8 \text{ mm}^2\text{C s}^2 & x_3 \leq 5\\ 4\text{e}6 \text{ mm}^2\text{C s}^2 & x_3 > 5 \end{cases},\] 
corresponding to mild carbon steel and its solid corrosion products, assumed to be iron (III) oxide, Fe$_2$O$_3$. We perform 50 PCG solutions with $\Delta t = 0.01$ seconds and a relative residual tolerance of $10^{-6}$. We vary the density of the tetrahedral mesh over the rectangular prism and measure the wall clock time per PCG iteration for each method. This provides a performance metric that accounts for computation time as well as the transfer of memory, and is consistent with the reporting of similar experiments \cite{Martinez-Frutos_CS, Martinez-Frutos_FEAD}. The results of these experiments are shown in Figure \ref{fig:PCG_time_per_iteration}. Although there is sublinear progression with coarse meshes as memory transfer time dominates, all of the implementations exhibit linear increase in computation with the number of degrees of freedom. The rate of increase of wall clock time for each GPU implementation, taken from the linear asymptote, are compared with the sparse CPU method in Table \ref{tab:Speedup}. We note that these factors may differ from the factors taken from comparisons using similar implementations in the literature  \cite{Martinez-Frutos_CS} for a number of reasons: lack of access to specific algorithmic details, hardware and API differences, and because our CPU baseline implementation is already highly optimized for this class of problem. 

We also include, for comparison, the same profiling results when the GPU methods are set to compute with single-precision floating point arithmetic. As seen in Figure \ref{fig:PCG_time_per_iteration}, the total timer per iteration decreases for all methods. We note that the point at which the dual GPU implementation becomes feasible is delayed to roughly double the problem size, as the balance between computation and memory transfer time is shifted. As the necessary computations become faster with single-precision arithmetic, the memory transfer overhead is more significant of a factor.

\begin{table}[h!]
\centering
\begin{tabular} {| c | c | c | c | c |}
\hline
Method & 1 & 2 & 3: single GPU & 3: dual GPU\\
\hline
Speedup Factor &  6.5 & 1.2 & 8.4 & 16.8\\
\hline
 \end{tabular}
\vspace{.2cm}
\captionsetup{margin = 0.5cm}
  \caption{Ratio of the rate of increase of wall clock time with increasing degrees of freedom between the sparse matrix CPU method and the linear asymptote of each GPU method.}
\vspace{-.4cm}
\label{tab:Speedup}
\end{table}

The experimentally observed linear scaling in computation with the number of degrees of freedom is expected based on the algorithms for computing MVM. Furthermore, it is seen that the efficiency of an implementation is directly related to the its care in the treatment of loading data from global memory. The second implementation has advantages over the first in that it does not require loading elementary matrix data or a second pass to sum contributions to a MVM. However, the cost is that it involves loading nine separate strips of data from the input vector to provide full information about all elements that are adjacent to a vertex. We see that this overhead outweighs the other merits of the implementation. The third implementation successfully adapts the memory management benefits of the previous two. A MVM requires only four strips of data from the input vector for the first pass, which are coalesced into four reads from global memory, and the contributions are partially summed so that the second pass will have less work. For the smallest problems considered here, the sparse matrix method is fastest, and the implementation on dual GPUs is the slowest. Further investigation into the dual GPUs implementation is done in the following section. 

\subsection{Multiple GPUs}
Performing the third PCG method on dual GPUs gives better performance for large problems, as expected. We see in Figure \ref{fig:PCG_time_per_iteration} that for problems with fewer than about 400,000 degrees of freedom, the increased memory transfer costs dominate total computation time. For larger problem, the burden of computing MVMs grows, so that increased parallelism is worth these costs. We use OpenCL profiling operations to confirm that the time spent actively computing MVMs is reduced by the same factor as the workload sharing (including the overlapping region), even though the directly observed total time per iteration exhibits smaller reduction. The effective doubling of performance reported in Table \ref{tab:Speedup} is taken from the relative slopes of the linear asymptote of each profiling curve. In other words, for the limit of problems in which computation dominates communication overhead. For problems larger than those considered here, it is expected that even more GPUs become practical, and we note that the extension is straightforward with the additive Schwartz method for domain decomposition.

\subsection{Further CPU Comparison}
\label{sec:serialComparison}
\begin{figure}
\centering
\includegraphics[width=\textwidth]{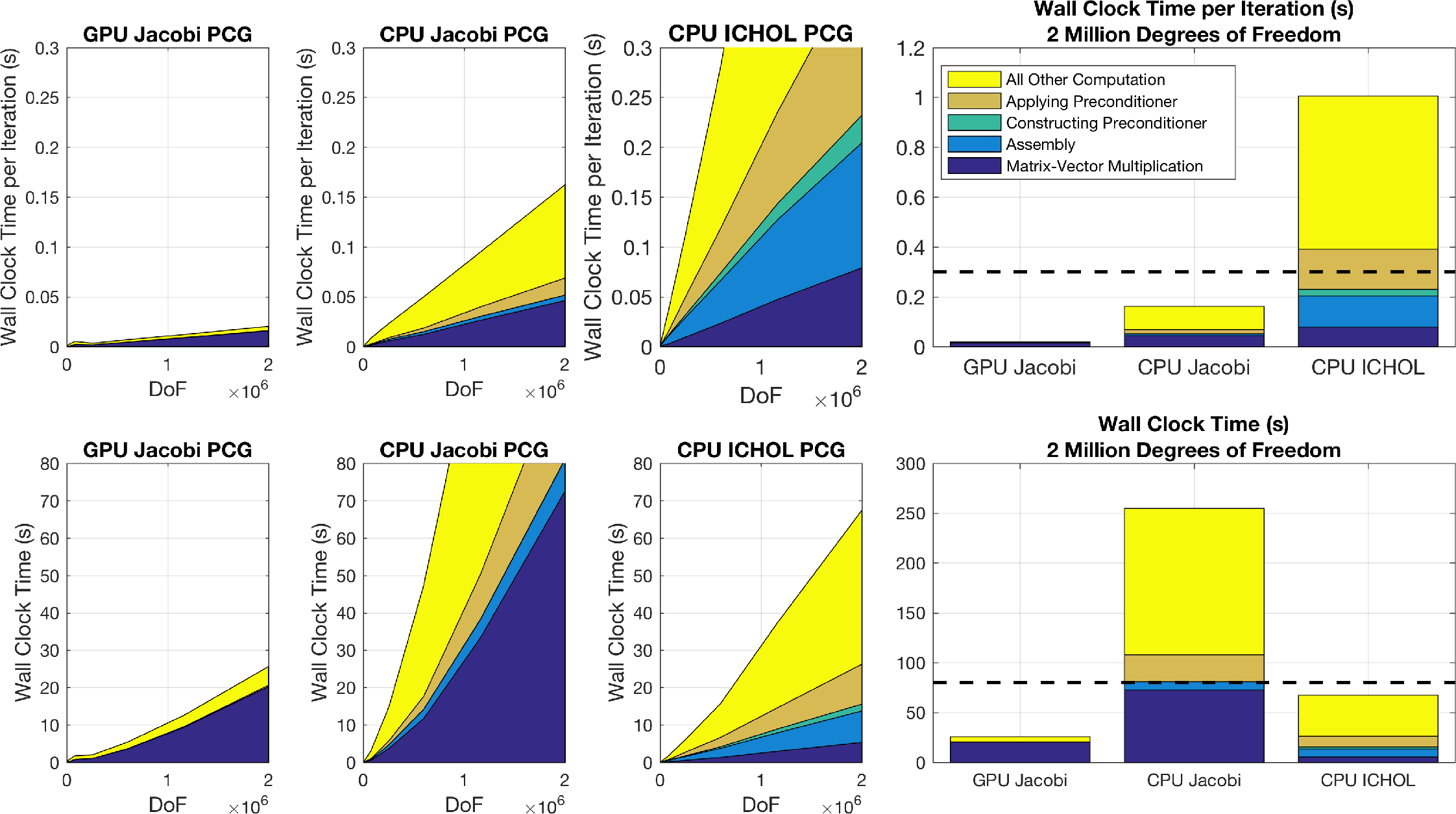}
\captionsetup{margin = 0.5cm}
\caption{Breakdown of major computational steps for 50 time steps using three algorithms. Time per iteration and total time are shown in the left panels across a range of mesh sizes. The right panels give a closer comparison of the three algorithms for a domain with 2 million degrees of freedom. The dashed black lines denote the maximum value of the time axis of the left plots.}
\vspace{-.4cm}
\label{fig:profiling_breakdown}
\end{figure}

We continue with a deeper comparison between our best performing heterogeneous computing method with a more sophisticated serial sparse matrix method. Since a standard CPU sparse matrix method can access the system matrix quickly, it is feasible to compute a better preconditioner matrix than the Jacobi preconditioner. This can take more time to produce and compute inverse matrix-vector multiplications with, but can vastly reduce the number of iterations required. We find that an incomplete Cholesky factorization with drop tolerance of $10^{-3}$ performs well for the simulations discussed in this paper \cite{golub2012matrix}. Table \ref{tab:ichol} shows the benefits of using this preconditioner for a transient heat conduction boundary value problem over 50 time steps. Since it is no longer fair to compare the average time per PCG iteration, we report the total time that both methods take to provide a solution, starting from knowledge of the elemental assembly matrices and a parameterization of material properties, including the incomplete Cholesky factorization for the serial method. The effective time per iteration is determined from this adjusted total time, rather than exclusively from the time in the PCG outer loop. A graphical breakdown of the computations involved is presented in Figure \ref{fig:profiling_breakdown}.

\begin{table}
\hspace{-2cm}
\begin{tabular} {| p{4cm} | c | c | c | c | c | c |}
\hline
Degrees of Freedom & 10.5$\times 10^{3}$  & 78.1$\times 10^{3}$  & 257$\times 10^{3}$  & 600$\times 10^{3}$  & 1.16$\times 10^{6}$  & 2.00$\times 10^{6}$ \\
\hline
\textbf{Sparse matrix CPU:} \newline Total time (s) & 0.2 & 1.6 & 6.2 & 16 & 37 & 67  \\
\hline
 Total Iterations & 39 & 47 & 54 & 56 & 62 & 67  \\
 \hline
Time per Iteration (s) & 5.6$\times 10^{-3}$ & 34$\times 10^{-3}$ & 110$\times 10^{-3}$ & 280$\times 10^{-3}$ & 600$\times 10^{-3}$ & 1000$\times 10^{-3}$  \\
\hline
\textbf{FG DbD:} \newline Total time (s) & 0.48 & 1.7 & 2.0 & 5.4 & 13 & 26 \\
\hline
 Total Iterations & 287 & 344 & 567 & 780 &1047  & 1278  \\
 \hline
Time per Iteration (s) & 1.6 $\times 10^{-3}$ & 5$\times 10^{-3}$ & 3.6$\times 10^{-3}$ & 6.9$\times 10^{-3}$ & 12$\times 10^{-3}$ & 20$\times 10^{-3}$  \\
\hline
 \end{tabular}
\vspace{.2cm}
\captionsetup{margin = 0.5cm}
  \caption{Full performance comparison between the sparse matrix CPU method with incomplete Cholesky preconditioning and our best performing implementation.}
\vspace{-.4cm}
\label{tab:ichol}
\end{table}

We see that the sparse matrix CPU method with incomplete Cholesky preconditioning takes similar overall time for modestly large systems. However, the differences in wall clock time diverge as explicit storage of the system matrix becomes more demanding. The heterogeneous computing method developed in this research outperform even sophisticated serial algorithms. 

\subsection{3D Coefficient Inverse Problem}
\label{sec:MCMC}
An important benefit of the currently proposed methods is that many simulations of heat conduction can be rapidly performed, in sequence, over a domain with varying thermal properties. This scenario arises in the solution of coefficient inverse problems, such as determining internal properties of a structure from noisy temperature measurements at its surface, in response to a known energy input. We explore such a case that has been earlier studied in two dimensions, motivated by a real-world corrosion detection problem \cite{EarlsMSSP}. By extending previous analysis methods from the literature to a full 3D model, we not only approach a more physically realistic use case, but are also able to relax certain symmetry restrictions on the boundary heat flux $f$. As a result, we found that an internal corrosion profile can be recovered with higher confidence, while using a heat source only one tenth as powerful as what was required in earlier work \cite{EarlsMSSP}. The details of the numerical experiment are now discussed.

\begin{wrapfigure}[26]{r}{0.5\textwidth}
\centering
\centerline{\includegraphics[width=0.5\textwidth]{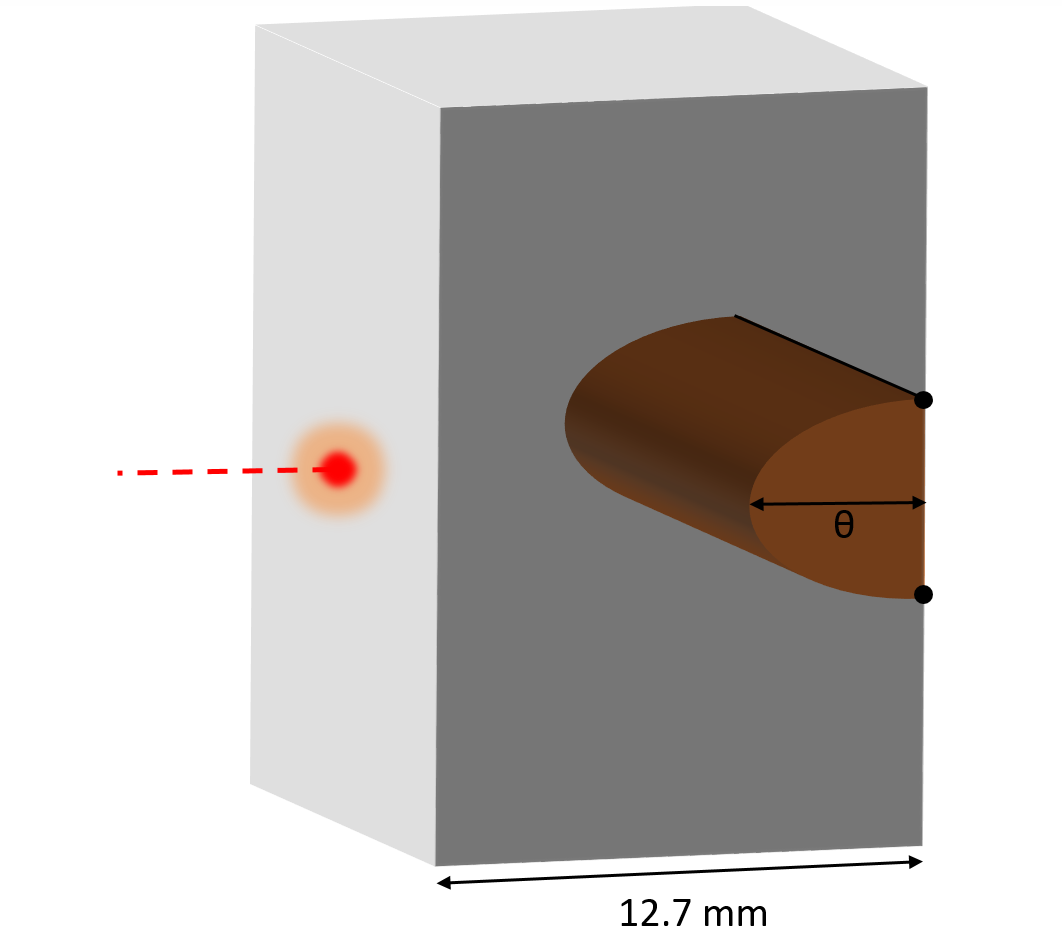}}
\captionsetup{margin = 0.1cm}
\caption{Parameterized corrosion pattern within a steel gusset plate. The parabolic corrosion boundary is ``anchored'' at the points shown in black, and grows out away from the rear boundary. Heat is input to the system on the front face, where the resulting temperature profile is also recorded.}
\label{fig:corrosion}
\end{wrapfigure}

Corrosion may form in a bridge structure in the crevice where a lower truss chord member meets a flat steel gusset plate connection element. It is observed that this corrosion within the gusset plate has a well defined geometric form, constant along the horizontal length of the truss chord, and varying as a quadratic function in the vertical direction \cite{EarlsMSSP}. The severity of the corrosion is parameterized by the penetration depth of the apex of the parabola into the steel, $\theta$, illustrated in Figure \ref{fig:corrosion}. We consider a case of corrosion that has penetrated 3.175 mm into a plate that is 12.7 mm thick: corresponding to a 25\% section loss. The relevant thermal properties are the same as in Section \ref{sec:performanceComparison}: A thermal input is furnished in the form of a 10 W laser beam having a Gaussian profile with 2 mm beam width. The structure is heated for $T_F = 10$ seconds, after which the surface temperature is recorded with an thermal camera. Surrogate field data, $\mathcal{D}$, are created with a high fidelity forward model (1.16 million degrees of freedom and 0.01 second time steps), interpolating the FEM solution to a finer rectangular grid, averaging the temperatures over each pixel area, and contaminating the result with noise; all so as to approximate a plausible digital thermal camera measuring device. For consistency with previous work, the camera is assumed to follow a noise model consisting of contaminating each pixel with independent and identically distributed Gaussian noise with zero mean and standard deviation of 0.1 $^\circ$C, then rounding the result to the nearest 0.1 $^\circ$C.

The inverse problem is solved with Markov chain Monte Carlo (MCMC), a Bayesian inference method \cite{Hastings}. The solution comes in the form of dependent samples taken from a probability distribution over the corrosion penetration depth, $p(\theta | \mathcal{D})$, known as the \textit{posterior} distribution. Since there is randomness in the surface temperature measurements, MCMC is able to automatically propagate a measure of uncertainty to the posterior distribution. The mean and standard deviation of a large number of posterior samples provide useful information about the underlying corrosion depth. The consequence is that every sample comes at the cost of performing a simulation with the candidate corrosion depth, $\hat{\theta}$, and evaluating the probability that the resulting FEM solution gave rise to the observed temperature data after passing through the camera noise model. This is known as the \textit{likelihood}, $p(\mathcal{D} | \hat{\theta})$. Assuming an ``uninformative'' uniform distribution over $\theta$, that is called a \textit{prior}, the corrosion depth is equally likely to be anywhere within the gusset plate thickness, Bayes' theorem states that the posterior distribution is proportional to the likelihood distribution, $p(\theta | \mathcal{D}) \propto p(\mathcal{D} | \theta)$, only scaled by a constant. Thus, the problem of estimating the posterior distribution can be reduced to sampling candidate values for $\theta$ and evaluating their corresponding likelihoods.

\begin{figure}
\centering
\includegraphics[width=\textwidth]{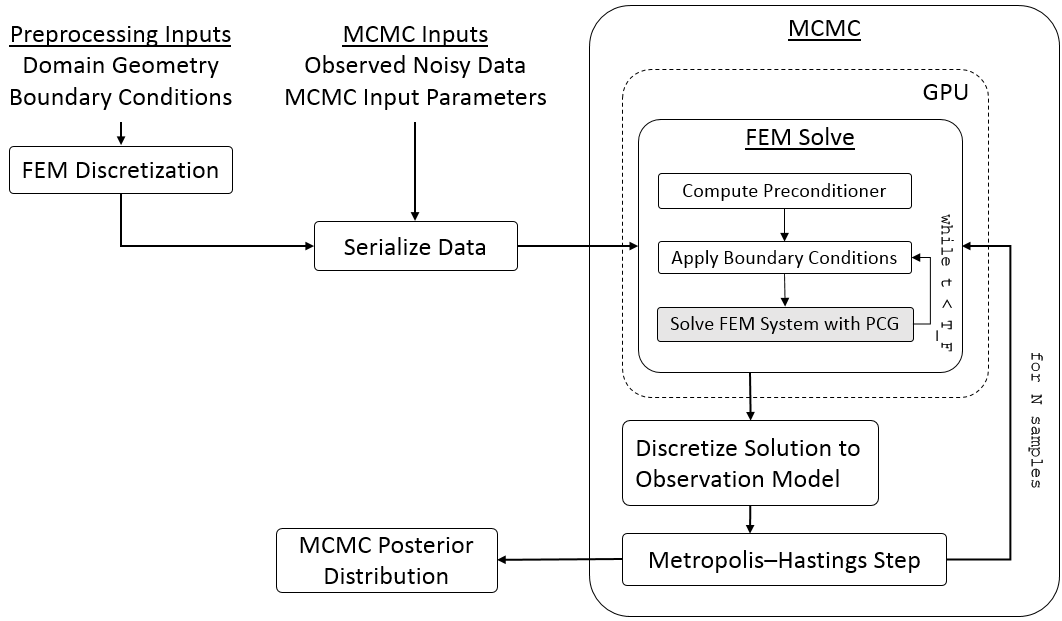}
\captionsetup{margin = 0.5cm}
\caption{Algorithmic flowchart for solving a coefficient inverse problem with MCMC and the heterogeneous computing FEM methods described in this work.}
\vspace{-.4cm}
\label{fig:flowchart}
\end{figure}

Markov chain Monte Carlo provides a systematic algorithm for producing the desired samples over $\theta$. Starting from an initial guess, $\theta_0$, a new candidate, $\hat{\theta}$, is randomly generated from some small neighborhood of $\theta_0$. The likelihoods for each value, given the observed data, are calculated through the use of the FEM solver. If the new sample is more likely to have produced the data, then it is set as $\theta_1$. Even if $\hat{\theta}$ was less likely to have produced the observed data, it may be randomly selected as the next sample with probability $p(\mathcal{D} | \hat{\theta})/p(\mathcal{D} | \theta_0)$. Otherwise, $\theta_0$ advances, repeated as sample $\theta_1$. This algorithm, known as Metropolis-Hastings sampling \cite{Hastings}, is summarized as 
\[ \theta_{i+1} = \begin{cases} \hat{\theta} &\text{ with probability } \min(1, p(\mathcal{D} | \hat{\theta})/p(\mathcal{D} | \theta_i)) \\ 
\theta_i & \text{ otherwise.} \end{cases} \]
Under some technical conditions that are easily verified for our system \cite{Link}, the sequence $\{\theta_i\}_{i = 1,2,3,...}$ is guaranteed to converge to samples from the desired distribution $p( \theta | \mathcal{D})$. The sampling process is summarized in Figure \ref{fig:flowchart}. The FEM implementation developed in this work is nested within two loops, as the critical stage of computation. Our GPU approach is massively parallelized for this task, and is designed to provide rapid successive solutions with only the transfer of $\theta_i$ from the host to the device (i.e. very low communication overhead). Furthermore, if the elemental assembly data is precomputed and loaded onto the device, the entire process can be carried out with minimal memory transfer.

\begin{figure}
\centering
\includegraphics[width=\textwidth]{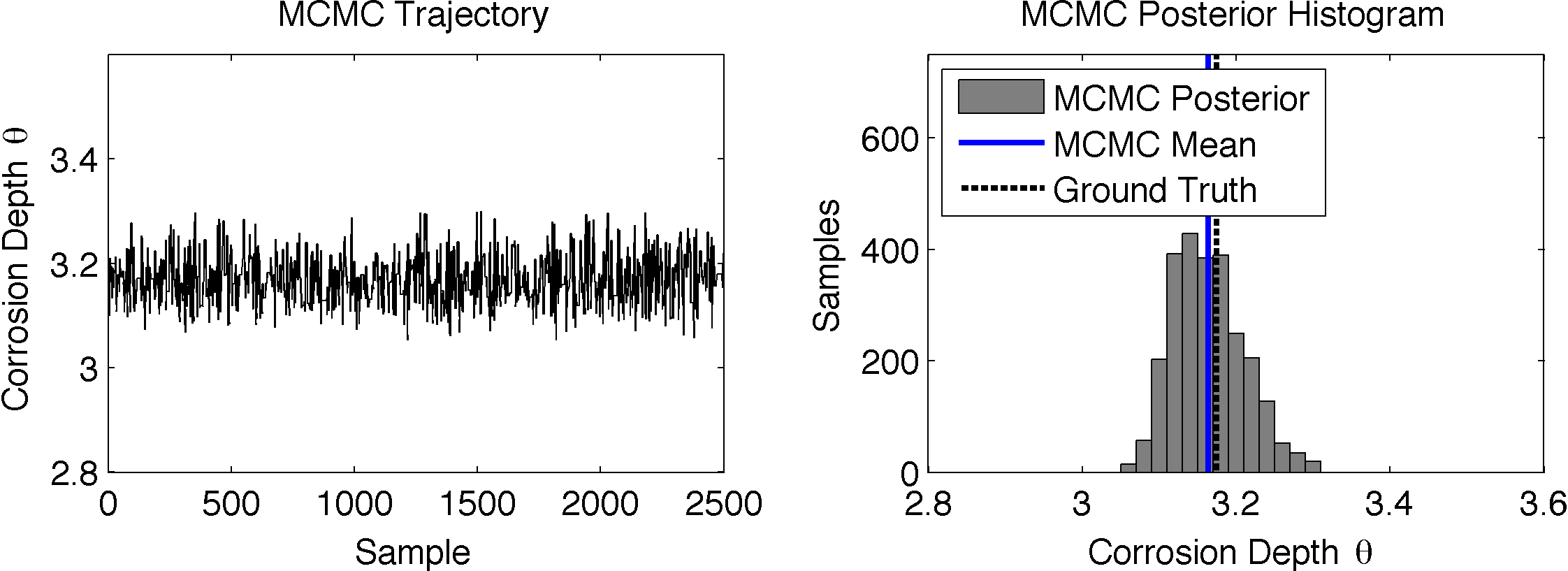}
\captionsetup{margin = 0.5cm}
\caption{(left) Trajectory of the 2500 sample Markov chain. (right) The same samples plotted as a histogram which approximates the distribution $p(\theta | \mathcal{D})$.}
\vspace{-.4cm}
\label{fig:MCMC}
\end{figure}

It is customary when using MCMC to discard some number of initial samples, to allow the chain to ``burn in'' to the posterior distribution, and away from its arbitrary initial value. For these experiments, we perform 200 burn-in iterations, starting from the middle of the prior distribution, and then record the subsequent 2500 samples. The trajectory of the Markov chain and a histogram of the samples are shown in Figure \ref{fig:MCMC}. It can be seen that the samples fall tightly around the ground truth value of 3.175 mm, denoting a confident solution to the coefficient inverse problem. The mean estimate is $3.16$ mm with standard deviation $0.05$ mm.

The demonstrated sequence of transient heat conduction simulations completed in 9 hours. If the same simulations are attempted using FEniCS, a modern open-source FEM computing platform \cite{FEniCSBbook}, the system runs out of memory in the assembly process. Computing on a series of smaller problems and extrapolating the run times (ignoring all contributions from just-in-time compilation) yields an expected total time of 142 hours. Using the in-house sparse matrix CPU code that was written specifically to solve this kind of sequence of heat conduction simulation, and discussed in Section \ref{sec:serialComparison}, would have taken approximately 38 hours. 

\section{Conclusions}
\label{sec:conclusions}
In this paper we described the development over three implementations of fine-grained assembly-free finite element simulation of heat transfer through heterogeneous media with the use of heterogeneous computing. The methods are motivated in terms of the interpretation of the FEM assembly operator as well as their particular implementation on GPUs. The implementation which gave the best performance was adapted for its simultaneous application on multiple GPUs. The ease of application of these methods to problems with spatially varying material properties is a direct consequence of the decomposition of the problem into local geometric components, rather than focusing on a global system matrix. Furthermore, successive simulations have low overhead, as mesh generation and computation of elemental assembly matrices can be done in advance. For these reasons, the methods that are developed in this research are particularly well-suited as forward models for coefficient inverse problems and topology optimization problems involving transient heat conduction. 

We finally note that the algorithm design methods that we have described are not restricted to the heat equation, and can be used as guidelines for solving other PDEs with spatially varying coefficients. To promote the use and extension of our work, we have used the non-proprietary OpenCL API and made the code for our best-performing implementation (both single and dual GPU versions) publicly available \cite{LoebGit}. 

\section{Acknowledgments}
This work was supported by the National Science Foundation (DGE-1144153).

\appendix
\section{GPU PCG Algorithm Details}
\label{app:PCG}
This appendix contains descriptions of the custom kernels that were written to perform linear algebra operations, how they use GPU memory, and how they fit into the PCG algorithm.

\subsection{Memory}
We first introduce the FEM data and indexing variables that are used throughout our kernels:
\begin{itemize}
\item M and F: arrays of data for elemental assembly matrices. Either the single matrices for fixed grid methods, or all entries for every element for the general methods
\item DoFMapLocal and coordinateMapLocal: small arrays of indices based on the geometry of the mesh to quickly determine local DoF to element assignment 
\item C: array of the number of divisions of the domain in each dimension
\item vert\_scale: array of the minimum vertex position and spacing between vertices in each dimension. Along with C, this allows the determination of the absolute location of a vertex given its global index. For a non-fixed grid approach, the user may start with a uniform mesh and deform it in a predetermined way (e.g. quadratically scaling the vertex locations), so that information can be used inside a kernel to still enable the recovery of absolute vertex position.
\item corr\_bound: spatial information about what part of the domain is corroded (or in general, a different material). For uniform corrosion after a certain depth, this can be that distance in the $z$-direction. For elliptical corrosion pits, it could contain the coordinates of the center of the ellipse and its axis lengths. The kernel must be programmed to know how to interpret this. 
\item mat\_coefs: values for $\rho C$ and $k$ for the different materials
\end{itemize}

Memory is allocated in GPU global memory buffers with flags to specify how the host and the device will access them. The flags are self-explanatory, and are combined with logical OR. The combinations used here are:
\begin{itemize}
\item HOST\_TO\_DEVICE\_COPY = (READ\_ONLY  $|$ HOST\_WRITE\_ONLY $|$ COPY\_HOST\_PTR)
\item HOST\_TO\_DEVICE\_USE  = (READ\_ONLY  $|$ HOST\_WRITE\_ONLY $|$ USE\_HOST\_PTR)
\item HOST\_READ\_WRITE     = (READ\_WRITE $|$ COPY\_HOST\_PTR)
\item PINNED              = (READ\_WRITE $|$ USE\_HOST\_PTR)
\item DEVICE\_READ\_WRITE   = (READ\_WRITE $|$ HOST\_NO\_ACCESS)
\end{itemize}

The memory buffers allocated for the FG DbD method are enumerated in Table \ref{tab:buffers}. Other methods do not differ much at this level. 

\begin{table}[h!]
\centering
\begin{tabular}{| c | c |  c | }
\hline
\textbf{Name} & \textbf{Flag} & \textbf{Initialization} \\ \hline
M\_FG\_buf & HOST\_TO\_DEVICE\_USE & M\_FG \\ \hline
K\_FG\_buf & HOST\_TO\_DEVICE\_USE & K\_FG\\ \hline
DoFMapLocal\_buf & HOST\_TO\_DEVICE\_COPY & DoFMapLocal\\ \hline
coordinateMapLocal\_buf & HOST\_TO\_DEVICE\_COPY & coordinateMapLocal\\ \hline
C\_buf & HOST\_TO\_DEVICE\_COPY & C\\ \hline
vert\_scale\_buf & HOST\_TO\_DEVICE\_COPY & vert\_scale\\ \hline
corr\_bounds\_buf & HOST\_TO\_DEVICE\_COPY & corr\_bounds\\ \hline
mat\_coefs\_buf & HOST\_TO\_DEVICE\_COPY & mat\_coefs\\ \hline
P\_buf & DEVICE\_READ\_WRITE & x.nbytes\\ \hline
u0\_buf & HOST\_TO\_DEVICE\_COPY & u0\\ \hline
Fdt\_buf & HOST\_TO\_DEVICE\_COPY & Fdt\\ \hline
x\_buf & HOST\_READ\_WRITE & x\\ \hline
b\_buf & DEVICE\_READ\_WRITE & x.nbytes\\ \hline
r\_buf & DEVICE\_READ\_WRITE & x.nbytes\\ \hline
d\_buf & DEVICE\_READ\_WRITE & x.nbytes\\ \hline
q\_buf & DEVICE\_READ\_WRITE & x.nbytes\\ \hline
s\_buf & DEVICE\_READ\_WRITE & x.nbytes\\ \hline
delta\_buf & PINNED & delta\\ \hline
alpha\_buf & DEVICE\_READ\_WRITE & delta.nbytes\\ \hline
neg\_alpha\_buf & DEVICE\_READ\_WRITE & delta.nbytes\\ \hline
delta\_new\_buf & DEVICE\_READ\_WRITE & delta.nbytes\\ \hline
beta\_buf & DEVICE\_READ\_WRITE & delta.nbytes\\ \hline
r1\_buf & DEVICE\_READ\_WRITE & r1\_size$\times$4\\ \hline
r2\_buf & DEVICE\_READ\_WRITE & r2\_size$\times$4\\ \hline
VVM\_loc\_buf & LocalMemory & max\_wg\_size$\times$4 \\ \hline
Ax\_split\_buf & DEVICE\_READ\_WRITE  & x.nbytes$\times$4\\ \hline
P\_split\_buf & DEVICE\_READ\_WRITE  & x.nbytes$\times$4\\ \hline
x\_local\_buf & LocalMemory & 32$\times$4$\times$4\\ \hline
\end{tabular}
\vspace{.2cm}
\captionsetup{margin = 0.5cm}
\caption{Memory buffer types and initialization values. If a number is given for initialization, the specified number of bytes is allocated.}
\vspace{-.4cm}
\label{tab:buffers}
\end{table}

\subsection{Kernels}
\subsubsection{VVM\_A}
First stage of a vector-vector multiplication. Takes pairwise scalar products of vector elements and then uses a standard parallel ``reduction'' algorithm to sum the results in $\mathcal{O}(\log(N))$ complexity. The work group size is maximized for the available hardware, and the partial sum  is reduced by a factor of two at each step. The total number of elements in the partial sum can only be reduced by a factor of the maximum work group size with a single kernel call. In case the length of the vector, nVertices, is not a multiple of the maximum work group size, max\_wg\_size, the number of work groups is rounded up. The size of the result (number of elements) is then 
\[ \mbox{r1\_size } =  \left\lceil \frac{\mbox{nVertices}}{\mbox{max\_workgroup\_size}} \right\rceil . \]

\begin{table}[h!]
\begin{tabular}{| c | c |  c | c | c |}
\hline
\textbf{Argument} & \textbf{In/Out} & \textbf{Memory space} & \textbf{Data type} & \textbf{Size}\\ \hline
*x & input & global & float64 & nVertices\\ \hline
*y & input & global & float64 & nVertices\\ \hline
nVertices & input &  & uint32 & 1\\ \hline
*VVM\_loc && local & float64 & max\_workgroup\_size\\ \hline
*r & output & global & float64 & r1\_size \\ \hline
\end{tabular}\\
\begin{tabular}{| c | c |  c | c |}
\hline
\textbf{global size} &r1\_size $\times$ max\_wg\_size\ & \textbf{local size} & max\_wg\_size\\ \hline
\end{tabular}
\end{table}

\subsubsection{VVM\_reduce}
Intermediate and final stages of vector-vector multiplication. Takes an intermediate partial sum from VVM\_A or itself and reduces it by a factor of max\_workgroup\_size. If nVertices $\leq$ max\_workgroup\_size, the vector-vector multiplication is completed. Otherwise further calls with this kernel are made. Since this can be iterative, set
\[ \mbox{rk\_size } = \left\lceil \frac{\mbox{r(k-1)\_size}}{\mbox{max\_workgroup\_size}} \right\rceil . \]
\begin{table}[h!]
\begin{tabular}{| c | c |  c | c | c |}
\hline
\textbf{Argument} & \textbf{In/Out} & \textbf{Memory space} & \textbf{Data type} & \textbf{Size}\\ \hline
*rPrevious & input & global & float64 & r(k-1)\_size $\times$ max\_workgroup\_size\\ \hline
nVertices & input &  & uint32 & 1\\ \hline
*VVM\_loc && local & float64 & max\_workgroup\_size\\ \hline
*r & output & global & float64 & rk\_size \\ \hline
\end{tabular}\\
\begin{tabular}{| c | c |  c | c |}
\hline
\textbf{global size} &rk\_size $\times$ max\_wg\_size\ & \textbf{local size} & max\_wg\_size\\ \hline
\end{tabular}
\end{table}

\subsubsection{VVM\_C}
An alternative final stage for vector-vector multiplication for step 5 of PCG. Rather than storing the final result of the sum, the scalar delta\_new is loaded and alpha = delta\_new/(d$^T$q) is stored, along with negative alpha.
\begin{table}[h!]
\begin{tabular}{| c | c |  c | c | c |}
\hline
\textbf{Argument} & \textbf{In/Out} & \textbf{Memory space} & \textbf{Data type} & \textbf{Size}\\ \hline
*rPrevious & input & global & float64 & r(k-1)\_size $\times$ max\_workgroup\_size\\ \hline
*delta & input & global & float64 & 1\\ \hline
nVertices & input &  & uint32 & 1\\ \hline
*b && local & float64 & max\_workgroup\_size\\ \hline
*alpha & output & global & float64 & 1\\ \hline
*neg\_alpha & output & global & float64 & 1\\ \hline
\end{tabular}\\
\begin{tabular}{| c | c |  c | c |}
\hline
\textbf{global size} &rk\_size\ & \textbf{local size} & max\_wg\_size\\ \hline
\end{tabular}
\end{table}

\subsubsection{VAVSP}
Computes the elementwise sum of a vector and a scalar multiplied by another vector. The scalars are loaded from GPU global memory, so a negative scalar must be used if vector subtraction is desired.
\begin{table}[h!]
\begin{tabular}{| c | c |  c | c | c |}
\hline
\textbf{Argument} & \textbf{In/Out} & \textbf{Memory space} & \textbf{Data type} & \textbf{Size}\\ \hline
*x & input & global & float64 & nVertices\\ \hline
*y & input & global & float64 & nVertices\\ \hline
*a &input & global & float64 & 1\\ \hline
*x\_plus\_ay & output & global & float64 & nVertices\\ \hline
\end{tabular}\\
\begin{tabular}{| c | c |  c | c |}
\hline
\textbf{global size} &nVertices & \textbf{local size} & None\\ \hline
\end{tabular}
\end{table}

\subsubsection{DIMVM}
Computes the matrix-vector multiplication where the matrix is the inverse of a diagonal matrix P. Element i of the result is x\_i/P\_\{i,i\}.
\begin{table}[h!]
\begin{tabular}{| c | c |  c | c | c |}
\hline
\textbf{Argument} & \textbf{In/Out} & \textbf{Memory space} & \textbf{Data type} & \textbf{Size}\\ \hline
*P & input & global & float64 & nVertices\\ \hline
*x & input & global & float64 & nVertices\\ \hline
*Pinvx & output & global & float64 & nVertices\\ \hline
\end{tabular}\\
\begin{tabular}{| c | c |  c | c |}
\hline
\textbf{global size} &nVertices & \textbf{local size} & None\\ \hline
\end{tabular}
\end{table}

\subsubsection{Beta\_update}
Computes the coefficient beta and updates delta in storage. Used to avoid data transfer between device and host for the performance of a small calculation.
\begin{table}[h!]
\begin{tabular}{| c | c |  c | c | c |}
\hline
\textbf{Argument} & \textbf{In/Out} & \textbf{Memory space} & \textbf{Data type} & \textbf{Size}\\ \hline
*delta\_new & input & global & float64 & 1\\ \hline
*delta\_old & both & global & float64 & 1\\ \hline
*beta & output & global & float64 & 1\\ \hline
\end{tabular}\\
\begin{tabular}{| c | c |  c | c |}
\hline
\textbf{global size} &1 & \textbf{local size} & 1\\ \hline
\end{tabular}
\end{table}

\subsubsection{u0\_update}
Computes an initial guess for the next time step, u\_0$^+$, based on a linear extrapolation from the initial guess of the current time step, u\_0 and the PCG solution of the current time step, u\_new,
\[ \mbox{ u\_0$^+$ = u\_new + (u\_new - u\_0).}\]
\begin{table}[h!]
\begin{tabular}{| c | c |  c | c | c |}
\hline
\textbf{Argument} & \textbf{In/Out} & \textbf{Memory space} & \textbf{Data type} & \textbf{Size}\\ \hline
*u0 & both & global & float64 & nVertices\\ \hline
*u\_new & input & global & float64 & nVertices\\ \hline
\end{tabular}\\
\begin{tabular}{| c | c |  c | c |}
\hline
\textbf{global size} &nVertices & \textbf{local size} & None\\ \hline
\end{tabular}
\end{table}

\subsubsection{FGDbDMVM\_A}
The implementation of matrix-vector multiplication (including cases for $\mathbf{A}\vec{x}$ and $a\mathbf{A}\vec{x} + \vec{b}$) and the determination of $P$ is dependent on the assembly perspective. The kernels for the third implementation (FG DbD  with memory coalescing) are discussed here.

The kernel computes contributions to a matrix-vector multiplication from each element according to coalesced memory access limitations. Data is loaded from the input vector as described in Section \ref{sec:FGEbE} and stored in local memory. Then each work item determines its global element index by first finding the cube in which it belongs: global\_id/6 - group\_id (global\_id is an unsigned integer, so integer division automatically rounds down), and then the tetrahedron within the cube: local\_id mod 6. 

The cube index also corresponds to the global vertex index for its first corner. This is used with the arrays C and vert\_scale to determine the vertex's absolute position. Based on this, the array coordinateMapLocal helps give the $x,y$, and $z$ positions of the other vertices of the tetrahedron, so that they can be compared with corr\_bounds to determine for each vertex whether it is in the corroded region of the domain or not. Material coefficients are taken from mat\_coefs and averaged over the element for both $\mathbf{M}_e$ and $\mathbf{K}_e$. 

Next, each element reads the data it needs about the input vector from local memory, referring to DoFMapLocal for the proper indices. This has been delayed as long as possible to hide the latency from the global memory load. Dot products are taken with the local assembly matrices, and the results are summed with the proper coefficients in another local memory array, Ax\_split\_local, which has 12 entries for every DoF, corresponding to the 12 possible contributing elements in the $\pm x$ directions. Finally, the work items that are responsible for loading and storing data sum over Ax\_split\_local for their DoF and write to global memory. The resulting memory buffer has 4 entries for every global DoF, one for every quadrant in the $y$-$z$ plane. These are all filled out in turn by further work groups.

The total number of work items needed for this kernel is found by determining the total number of elements that need to be considered including the padding in $+x$ and $+y$ directions, dividing by 30 since every work group yields the elemental contributions for blocks of 30 vertices, and multiplying by the work group size. As with vector-vector multiplications, we round up the integer division
\[ \mbox{MVM\_global\_size} = 186 \left\lceil \frac{6(C[0] + 1)(C[1] + 1)C[2]}{180} \right\rceil . \]
\begin{table}[h!]
\begin{tabular}{| c | c |  c | c | c |}
\hline
\textbf{Argument} & \textbf{In/Out} & \textbf{Memory space} & \textbf{Data type} & \textbf{Size}\\ \hline
*M\_FG & input & constant & float64 & 4\\ \hline
*K\_FG & input & constant & float64 & 4\\ \hline
*x & input & global & float64 & nVertices\\ \hline
*DoFMapLocal & input & constnat & uint32 & 12 \\ \hline
*coordinateMapLocal & input & constant & float64 & 12 \\ \hline
*C & input & constant & uint32 & 3 \\ \hline
*vert\_scale & input & constant & float64 & 6 \\ \hline
*corr\_bounds & input & constant & float64 &  \\ \hline
*mat\_coefs & input & constant & float64 & 4 \\ \hline
theta & input &  & float64 & 1 \\ \hline
dt & input &  & float64 & 1 \\ \hline
*x\_local &  & local & float64 & 4$\times$32 \\ \hline
*Ax\_split\_local &  & local & float64 & 4$\times$32$\times$12 \\ \hline
*Ax\_split & output & global & float64 & nVertices$\times$4\\ \hline
\end{tabular}\\
\begin{tabular}{| c | c |  c | c |}
\hline
\textbf{global size} &MVM\_global\_size & \textbf{local size} & 186\\ \hline
\end{tabular}
\end{table}

\subsubsection{FGDbDMVM\_B}
Finishes the matrix-vector multiplication started by FGDbDMVM\_A. Work items sum the four contributions to each DoF from Ax\_split and store them in the final result array. 
\begin{table}[h!]
\begin{tabular}{| c | c |  c | c | c |}
\hline
\textbf{Argument} & \textbf{In/Out} & \textbf{Memory space} & \textbf{Data type} & \textbf{Size}\\ \hline
*Ax\_split & input & global & float64 & nVertices$\times$4\\ \hline
*Ax & output & global & float64 & nVertices\\ \hline
\end{tabular}\\
\begin{tabular}{| c | c |  c | c |}
\hline
\textbf{global size} &MVM\_global\_size & \textbf{local size} & 186\\ \hline
\end{tabular}
\end{table}

\subsubsection{FGDbDMVM\_C}
Finishes the matrix-vector multiplication started by FGDbDMVM\_A with an extra SAXPY operation so that a separate call to AXPY is not necessary. Work items sum the four contributions to each DoF from Ax\_split, multiply them by a scalar, add them to an element from another vector, and store the result. 
\begin{table}[h!]
\begin{tabular}{| c | c |  c | c | c |}
\hline
\textbf{Argument} & \textbf{In/Out} & \textbf{Memory space} & \textbf{Data type} & \textbf{Size}\\ \hline
*Ax\_split & input & global & float64 & nVertices$\times$4\\ \hline
*b & input & global & float64 & nVertices\\ \hline
c & input &  & float64 & 1\\ \hline
*cAx\_plus\_b & output & global & float64 & nVertices\\ \hline
\end{tabular}\\
\begin{tabular}{| c | c |  c | c |}
\hline
\textbf{global size} &MVM\_global\_size & \textbf{local size} & 186\\ \hline
\end{tabular}
\end{table}

\subsubsection{Jacobi\_A}
\label{app:Jacobi}
Determines contributions to the Jacobi preconditioner $P$. The algorithm is the same as FGDbDMVM\_A, except instead of taking dot products with elemental assembly matrices and an input vector, the diagonal elements of the elemental assembly matices are scaled according to material coefficients, combined, and stored. Calling FGDbDMVM\_B on the result produces the diagonal of $P$.\\ 
Arguments for this kernel are the same as for FGDbDMVM\_A, except without the need for *x and *x\_local.

\subsection{PCG Again}
The kernels and memory usage is as follows, folowing the steps in Section \ref{PCG}.  The syntax below is a small abbreviation of the actual PyOpenCL code, and has the form \\
kernel\_instance = kernel\_name(args)

First compute $P$ and the right hand side vector $\vec{b} = \mathbf{L}\vec{u}_i + \vec{F}$. 

\begin{itemize}
\item knl\_PA = Jacobi\_A(M\_FG\_buf, K\_FG\_buf, DoFMapLocal\_buf, 
          coordinateMapLocal\_buf, C\_buf, vert\_scale\_buf, corr\_bounds\_buf, mat\_coefs\_buf,
          theta,  dt, Ax\_split\_local\_buf, P\_split\_buf)

\item knl\_PB = FGDbDMVM\_B(P\_split\_buf, P\_buf)

\item knl\_RHS\_A = FGDbDMVM\_A(M\_FG\_buf, K\_FG\_buf, u0\_buf, DoFMapLocal\_buf, 
         coordinateMapLocal\_buf, C\_buf, vert\_scale\_buf, corr\_bounds\_buf, mat\_coefs\_buf, 
         (1-theta), dt, x\_local\_buf, Ax\_split\_local\_buf, Ax\_split\_buf)
\item knl\_RHS\_B = FGDbDMVM\_C(Ax\_split\_buf, Fdt\_buf, np.float64(1), b\_buf)
\end{itemize}
Setup for PCG
\begin{itemize}
\item knl\_1A = FGDbDMVM\_A(M\_FG\_buf, K\_FG\_buf, x\_buf, DoFMapLocal\_buf, 
         coordinateMapLocal\_buf, C\_buf, vert\_scale\_buf, corr\_bounds\_buf, 
         mat\_coefs\_buf, theta, dt, x\_local\_buf, Ax\_split\_local\_buf, Ax\_split\_buf)
\item knl\_1B = FGDbDMVM\_C(Ax\_split\_buf, b\_buf, np.float64(-1), r\_buf)

\item knl\_2 = DIMVM(P\_buf, r\_buf, d\_buf)

\item knl\_3A = VVM\_A(r\_buf, d\_buf, nVertices, VVM\_loc\_buf, r1\_buf)
\item knl\_3B = VVM\_reduce(r1\_buf, r1\_size, VVM\_loc\_buf, r2\_buf)
\item knl\_3C = VVM\_reduce(r2\_buf, r2\_size, VVM\_loc\_buf, delta\_buf)
\end{itemize}

Perform one PCG iteration
\begin{itemize}
\item knl\_4A = FGDbDMVM\_A(M\_FG\_buf, K\_FG\_buf, d\_buf, DoFMapLocal\_buf, 
         coordinateMapLocal\_buf, C\_buf, vert\_scale\_buf, corr\_bounds\_buf, 
         mat\_coefs\_buf, theta, dt, x\_local\_buf, Ax\_split\_local\_buf, Ax\_split\_buf)
\item knl\_4B = FGDbDMVM\_B(Ax\_split\_buf, q\_buf)

\item knl\_5A = VVM\_A(d\_buf, q\_buf, nVertices, VVM\_loc\_buf, r1\_buf)
\item knl\_5B = VVM\_reduce(r1\_buf, r1\_size, VVM\_loc\_buf, r2\_buf)
\item knl\_5C = VVM\_C(r2\_buf, delta\_buf, r2\_size, VVM\_loc\_buf, alpha\_buf, neg\_alpha\_buf)

\item knl\_6 = AXPY(x\_buf, d\_buf, alpha\_buf, x\_buf)

\item knl\_7A = FGDbDMVM\_A(M\_FG\_buf, K\_FG\_buf, x\_buf, DoFMapLocal\_buf, 
         coordinateMapLocal\_buf, C\_buf, vert\_scale\_buf, corr\_bounds\_buf, 
         mat\_coefs\_buf, theta, dt, x\_local\_buf, Ax\_split\_local\_buf, Ax\_split\_buf)
\item knl\_7B = FGDbDMVM\_C(Ax\_split\_buf, b\_buf, np.float64(-1), r\_buf)

\item knl\_7 = AXPY(r\_buf, q\_buf, neg\_alpha\_buf, r\_buf)

\item knl\_8 = DIMVM(P\_buf, r\_buf, s\_buf)

\item knl\_9A = VVM\_A(r\_buf, s\_buf, nVertices, VVM\_loc\_buf, r1\_buf)
\item knl\_9B = VVM\_reduce(r1\_buf, r1\_size, VVM\_loc\_buf, r2\_buf)
\item knl\_9C = VVM\_reduce(r2\_buf, r2\_size, VVM\_loc\_buf, delta\_new\_buf)

\item knl\_10 = Beta\_update(delta\_new\_buf, delta\_buf, beta\_buf)

\item knl\_11 = AXPY(s\_buf, d\_buf, beta\_buf, d\_buf)

\item knl\_u0 = u0\_update(u0\_buf, x\_buf)
\end{itemize}

The CPU handles logical decisions and calls these kernels according to the algorithm until convergence is realized.

\bibliography{References} 
\bibliographystyle{plain}

\end{document}